\documentclass{amsart}
\usepackage{amsmath, amssymb, latexsym, amscd, amsthm,amsfonts,amstext}
\usepackage[mathscr]{eucal}
\usepackage{graphicx}
\usepackage{subfig}

\theoremstyle{definition}

\theoremstyle{remark}

\numberwithin{equation}{section}

\input{tcilatex}

\begin{document}
\title[Convexification for a CIP for the RTE]{Convexification Numerical
Method for a Coefficient Inverse Problem for the Radiative Transport Equation%
}
\author{Michael V. Klibanov}
\address{Department of Mathematics and Statistics, University of North
Carolina at Charlotte, Charlotte, NC, 28223, USA}
\email{mklibanv@uncc.edu}
\author{Jingzhi Li}
\address{Department of Mathematics \& National Center for Applied
Mathematics Shenzhen \& SUSTech International Center for Mathematics,
Southern University of Science and Technology, Shenzhen 518055, P.~R.~China}
\email{li.jz@sustech.edu.cn}
\thanks{MK and LN were supported by US Army Research Laboratory and US Army
Research Office Grant W911NF-19-1-0044. LN was supported in part by NSF
Grant \#DMS-2208159 and by funds provided by the Faculty Research Grant
program at UNC Charlotte Fund \# 111272. JL was partially supported by the
NSF of China No. 11971221, Guangdong NSF Major Fund No. 2021ZDZX1001, the
Shenzhen Sci-Tech Fund No. RCJC20200714114556020, JCYJ20200109115422828 and
JCYJ20190809150413261.}
\author{Loc H. Nguyen}
\address{Department of Mathematics and Statistics, University of North
Carolina at Charlotte, Charlotte, NC, 28223, USA}
\email{loc.nguyen@uncc.edu}
\author{Zhipeng Yang}
\address{Department of Mathematics, Southern University of Science and
Technology, Shenzhen 518055, P.~R.~China}
\email{yangzp@sustech.edu.cn}
\subjclass[2020]{Primary 35R30; Secondary 65M32}
\keywords{coefficient inverse problem, radiative transport equation,
convexification method, a special orthonormal basis, global convergence
analysis, numerical experiments}

\begin{abstract}
An $\left( n+1\right) -$D coefficient inverse problem for the radiative
stationary transport equation is considered for the first time. A globally
convergent so-called convexification numerical \ method is developed and its
convergence analysis is provided. The analysis is based on a Carleman
estimate. In particular, convergence analysis implies a certain uniqueness
theorem. Extensive numerical studies in the 2-D case are presented.
\end{abstract}

\maketitle












\section{Introduction}

\label{sec:1}

The stationary radiative transfer equation (RTE) governs the propagation of
the radiation field in media with absorbing, emitting and scattering
radiation. RTE has broad applications in optics, including diffuse optical
tomography \cite{Heino}, astrophysics, atmospheric science and other applied
disciplines. For example, in the single particle emission tomography the
coefficient we reconstruct is the emission coefficient \cite[formula (2.1)]%
{Nat}. The RTE we consider here is a more general one since we introduce an
integral operator in it.

For the first time, we develop in this paper a globally convergent numerical
method for a Coefficient Inverse Problem (CIP) of the recovery of a
spatially distributed coefficient of RTE in the $\left( n+1\right) -$D case, 
$n\geq 1$. All CIPs are both nonlinear and ill-posed. Numerical methods for
inverse source problems for RTE were developed in \cite{HT1,HT2,HT3,SKN}. In
the case of single particle emission tomography, i.e. when the kernel of the
integral operator in RTE is the identical zero, an inversion formula for the
inverse source problem was first derived in \cite{Nov} and tested
numerically in \cite{GN}. Inverse problems of \cite{HT1,HT2,HT3,GN,Nov,SKN}
are linear ones.

Our CIP is formally determined, i.e. the number of free variables in the
data equals the number of free variables in the unknown coefficient. Our
data are incomplete, i.e. the source runs along an interval of a straight
line and the data are measured only at a part of the boundary $\partial
\Omega $ of the domain of interest $\Omega $. This is unlike the classical
case of X-ray tomography when the source runs all around $\Omega $ and the
data are measured on the whole boundary $\partial \Omega $.

A significant attention has been paid to the questions of uniqueness and
stability of CIPs for RTE, see, e.g. \cite{Bal,GY,KP,Lay,SU} and the
references cited therein. We refer here only to a limited number of works on
this topic since our interest in this paper is in a numerical development.

Our goal is to construct a globally convergent numerical method for our CIP,
to conduct its convergence analysis and to confirm our method via numerical
experiments. To achieve this goal, we develop here a new version of the
so-called $\emph{convexification}$ method \cite{KL}. The convexification
constructs a least squares cost functional, which is strictly convex on a
convex set in an appropriate Hilbert space. The diameter of this set is
fixed and is an arbitrary one. We prove existence and uniqueness of the
minimizer of that functional on that set and estimate convergence rate of
minimizers to the true solution depending on the level of the noise in the
data. As a by-product, we obtain a certain uniqueness result for our CIP,
see Theorem 4.3 below. Also, we establish the global convergence of the
gradient descent method of the minimization of our functional. Recall that
only local convergence of this method can be proven in the non-convex case.
In addition, we present numerical experiments in the 2-D case.

We call a numerical method for a CIP\ \emph{globally convergent} if a
theorem is proven, which claims that this method delivers at least one point
in a sufficiently small neighborhood of the solution of this CIP\ without an
advanced knowledge of that neighborhood, also, see \cite[Definiion 1.4.2]{KL}
for a similar statement. In other words, a good approximation for the true
solution can be obtained without a knowledge of a good first guess about
this solution.

Conventional numerical methods for CIPs are based on the minimizations of
least squares cost functionals, see, e.g. \cite{Chavent, Gonch2}. However,
these functionals are usually non-convex. This means, in turn that they
typically have multiple local minima and ravines. On the other hand, since
any gradient-like optimization method can stop at any point of a local
minimum, then this technique needs to have a good first guess about the
solution. Besides, there is no rigorous guarantee of neither the existence
of the global minimizer nor that such a \ minimizer, if it exists, is indeed
close to the true solution. These considerations have prompted the first
author to work on the developments of the convexification technique with the
first publications \cite{KlibIous,Klib97}.

The key element of the convexification is the presence of the so-called
Carleman Weight Function (CWF) in the above mentioned functional. The
presence of the CWF ensures the strict convexity property of that functional
on that bounded set. The CWF is the function, which is used as the weight in
the Carleman estimate for the corresponding differential operator, see, e.g.
books \cite{BY,KL} for Carleman estimates. In particular, we demonstrate in
our numerical experiments that the absence of the CWF leads to a significant
deterioration of the numerical solution.

The convexification uses the idea of the so-called Bukhgeim-Klibanov method
(BK), which is based on Carleman estimates. BK was originally introduced in
the field of CIPs in 1981 \cite{BukhKlib} only for the proofs of uniqueness
theorems for multidimensional CIPs. The work \cite{BukhKlib} has generated
many publications of many authors since then, see, e.g. books \cite{BY,KL},
papers \cite{GY,KP,Ksurvey,Lay} for a few of such publications as well as
references cited therein. The majority of currently known publications about
BK is also dedicated to the issues of uniqueness and stability of CIPs.
Unlike these, it was originally proposed in \cite{KlibIous,Klib97} to use
the idea of BK for the construction of the above mentioned globally strictly
convex cost functionals for CIPs.

While initially only analytical results for the convexification were
derived, an active exploration of this method for numerical studies of CIPs
has started from the work \cite{Bak}, which has removed some obstacles for
real computations, see, e.g., publications \cite{Khoa2,KLZ18,KL} for a few
examples out of many. In particular, \cite{Khoa2} and \cite[Chapter 10]{KL}
demonstrate the performance of the convexification for experimentally
collected data for backscatter microwaves. An updated version of the
convexification, which combines Carleman estimates with the contraction
principle, can be found in \cite{Baud1,Baud2,Buhan}.

In our approach, we use a certain approximate mathematical model. Our model
amounts to the truncation of a Fourier-like series with respect to a special
orthonormal basis in $L_{2}\left( a,b\right) $ as well as to the so-called
\textquotedblleft partial finite differences", see subsection 3.3 for
details. This basis was originally introduced in \cite{Klibanov:jiip2017},
also see \cite[Chapter 6]{KL}. We do not know how to prove convergence of
our method in the case when the number of terms of this series $N\rightarrow
\infty .$ Therefore, the \emph{only way} to verify the validity of this
model is via numerical simulations. In this regard, we note that truncations
of Fourier-like series with respect to the same basis were done for a
variety of CIPs in \cite{Khoa2,KLZ18,SKN} and in \cite[Chapters 7,10,12]{KL}%
. In each of these cases, the validity of the corresponding approximate
mathematical model was verified numerically well. Similar cases of truncated
Fourier series without proofs of convergence at $N\rightarrow \infty $ can
be observed for inverse problems considered by some other authors, see, e.g. 
\cite{GN,Kab}. And successful numerical verifications also took place in
these references.

Philosophically, the situation of our approximate mathematical model being
verified in numerical studies is somewhat similar with the well known
situation of the Huygens-Fresnel diffraction theory in optics. This theory
is not yet derived rigorously from the Maxwell's system. This is why it is
stated in section 8.1 of the classical textbook of Born and Wolf \cite{BW}
that \textquotedblleft \emph{because of mathematical difficulties,
approximate mathematical models must be used in most cases of practical
interest. Of these the theory of Huygens and Fresnel is by far most powerful
and adequate for the treatment of the majority of problems encountered in
experimental optics.}" The conclusion we draw from this is that in the case
of significant challenges for a certain applied mathematical problem, it is
reasonable to introduce an approximate mathematical model. However, this
model must be verified numerically.

All functions considered below are real valued ones. Let $B$ be a Banach
space of functions and $s\geq 2$ be an integer. Below $B_{s}=\underbrace{%
B\times B\times ...\times B},$ $s$ terms, is the space of $s-$dimensional
vector functions generated by $B$ and with the obvious extension of the $B-$%
norm.

In section 2 we pose forward and inverse problems and prove existence and
uniqueness theorem for the solution of the forward problem. In section 3 we
describe our transformation procedure. In section 4 we introduce our
convexification functional and formulate five theorems. These theorems are
proven in section 5. In section 6 we present our numerical results.

\section{Statements of Forward and Inverse Problems}

\label{sec:2}

For $n\geq 1,$ points in $\mathbb{R}^{n+1}$ are denoted below as $\mathbf{x}%
=(x_{1},x_{2},...,x_{n-1},y)\in \mathbb{R}^{n+1}.$ Let numbers $A,a,b,d>0$,
where 
\begin{equation}
1<a<b.  \label{2.0}
\end{equation}%
Define the rectangular prism $\Omega \subset \mathbb{R}^{n+1}$ and parts $%
\partial _{1}\Omega ,\partial _{2}\Omega ,\partial _{3}\Omega $ of its
boundary $\partial \Omega $ as: 
\begin{equation}
\Omega =\{\mathbf{x}:-A<x_{1},...,x_{n}<A,a<y<b\},  \label{2.1}
\end{equation}%
\begin{equation}
\partial _{1}\Omega =\left\{ \mathbf{x}:-A<x_{1},...,x_{n}<A,y=a\right\} ,
\label{2.2}
\end{equation}%
\begin{equation}
\partial _{2}\Omega =\left\{ \mathbf{x:}-A<x_{1},...,x_{n}<A,y=b\right\} ,
\label{2.3}
\end{equation}%
\begin{equation}
\partial _{3}\Omega =\left\{ x_{i}=\pm A,y\in \left( a,b\right)
,i=1,...,n,\right\} .  \label{2.4}
\end{equation}%
Let $\Gamma _{d}$ be the line where the external sources are, 
\begin{equation}
\Gamma _{d}=\{\mathbf{x}_{\alpha }=(\alpha ,0,...,0):\alpha \in \lbrack
-d,d]\}.  \label{2.40}
\end{equation}%
Hence, $\Gamma _{d}$ is a part of the $x_{1}-$axis. It follows from (\ref%
{2.0}) and (\ref{2.1}) that $\Gamma _{d}\cap \overline{\Omega }=\varnothing
. $

Let the points of external sources $\mathbf{x}_{\alpha }$ run along $\Gamma
_{d},$ $\mathbf{x}_{\alpha }\in \Gamma _{d}$. Let $\varepsilon >0$ be a
sufficiently small number. To avoid dealing with singularities, we model the 
$\delta \left( \mathbf{x}\right) -$function as:%
\begin{equation}
f\left( \mathbf{x}\right) =C_{\varepsilon }\left\{ 
\begin{array}{c}
\exp \left( \frac{\left\vert \mathbf{x}\right\vert ^{2}}{\varepsilon
^{2}-\left\vert \mathbf{x}\right\vert ^{2}}\right) ,\left\vert \mathbf{x}%
\right\vert <\varepsilon , \\ 
0,\left\vert \mathbf{x}\right\vert \geq \varepsilon%
\end{array}%
\right. ,  \label{2.5}
\end{equation}%
where the constant $C_{\varepsilon }$ is chosen such that%
\begin{equation}
C_{\varepsilon }\mathop{\displaystyle \int}\limits_{\left\vert \mathbf{x}%
\right\vert <\varepsilon }\exp \left( \frac{\left\vert \mathbf{x}\right\vert
^{2}}{\varepsilon ^{2}-\left\vert \mathbf{x}\right\vert ^{2}}\right) d%
\mathbf{x}=1.  \label{2.50}
\end{equation}%
Hence, the function $f\left( \mathbf{x}-\mathbf{x}_{\alpha }\right) =f\left(
x_{1}-\alpha ,x_{2},...,x_{n},y\right) \in C^{\infty }\left( \mathbb{R}%
^{n+1}\right) $ plays the role of the source function for the source $%
\mathbf{x}_{\alpha }.$ We choose $\varepsilon $ so small that 
\begin{equation}
f\left( \mathbf{x}-\mathbf{x}_{\alpha }\right) =0,\forall \mathbf{x}\in 
\overline{\Omega },\forall \mathbf{x}_{\alpha }\in \Gamma _{d}.  \label{2.6}
\end{equation}

Let $u(\mathbf{x},\alpha )$ denotes the steady-state radiance at the point $%
\mathbf{x}$ generated by the source function $f\left( \mathbf{x}-\mathbf{x}%
_{\alpha }\right) .$ We assume that the function $u(\mathbf{x},\alpha )$ is
governed by the RTE of the following form \cite{Heino}:%
\begin{equation*}
\nu (\mathbf{x,}\alpha )\cdot \nabla _{\mathbf{x}}u(\mathbf{x},\alpha
)+a\left( \mathbf{x}\right) u(\mathbf{x},\alpha )
\end{equation*}%
\begin{equation}
=\mu _{s}(\mathbf{x})\int_{\Gamma _{d}}K(\mathbf{x},\alpha ,\beta )u(\mathbf{%
x},\beta )d\beta +f\left( \mathbf{x}-\mathbf{x}_{\alpha }\right) ,\text{ }%
\mathbf{x}\in \mathbb{R}^{n+1},\mathbf{x}_{\alpha }\in \Gamma _{d}.
\label{2.7}
\end{equation}%
The kernel $K(\mathbf{x},\alpha ,\beta )$ of the integral operator in (\ref%
{2.7}) is called the \textquotedblleft phase function", 
\begin{equation}
K(\mathbf{x},\alpha ,\beta )\geq 0,\mathbf{x\in }\overline{\Omega };\text{ }%
\alpha ,\beta \in \left[ -d,d\right] ,\text{ }  \label{2.8}
\end{equation}%
see \cite{Heino}. In addition, we assume that 
\begin{equation}
K(\mathbf{x},\alpha ,\beta )\in C^{2}\left( \overline{\Omega }\times \left[
-d,d\right] ^{2}\right) .  \label{2.9}
\end{equation}%
In equation (\ref{2.7}),%
\begin{equation}
a\left( \mathbf{x}\right) =\mu _{a}\left( \mathbf{x}\right) +\mu _{s}(%
\mathbf{x}),  \label{2.10}
\end{equation}%
where $\mu _{a}\left( \mathbf{x}\right) $ and $\mu _{s}(\mathbf{x})$ are the
absorption and scattering coefficients respectively. As stated in section 1, 
$a\left( \mathbf{x}\right) $ is the emission coefficient \cite[formula (2.1)]%
{Nat}. We assume that 
\begin{equation}
\mu _{a}\left( \mathbf{x}\right) ,\mu _{s}(\mathbf{x})\geq 0,\mu _{a}\left( 
\mathbf{x}\right) =\mu _{s}(\mathbf{x})=0\text{, }\mathbf{x}\in \mathbb{%
\mathbb{R}}^{n+1}\setminus \Omega ,  \label{2.11}
\end{equation}%
\begin{equation}
\mu _{a}\left( \mathbf{x}\right) ,\mu _{s}(\mathbf{x})\in C^{2}\left( 
\mathbb{\mathbb{R}}^{n+1}\right) .  \label{2.12}
\end{equation}

For two arbitrary points $\mathbf{x},\mathbf{z}\in \mathbb{R}^{n+1}$ let $%
L\left( \mathbf{x},\mathbf{z}\right) $ be the line segment connecting these
points and let $ds$ be the element of the euclidean length on $L\left( 
\mathbf{x},\mathbf{z}\right) .$ In (\ref{2.7}) $\nu (\mathbf{x,}\alpha )$
denotes the unit vector, which is parallel to $L\left( \mathbf{x},\mathbf{x}%
_{\alpha }\right) ,$ 
\begin{equation}
\nu (\mathbf{x,}\alpha )=\frac{\mathbf{x}-\mathbf{x}_{\alpha }}{\left\vert 
\mathbf{x}-\mathbf{x}_{\alpha }\right\vert }=\frac{\left( x_{1}-\alpha
,x_{2},...,x_{n},y\right) }{\sqrt{\left( x_{1}-\alpha \right)
^{2}+x_{2}^{2}+...+x_{n}^{2}+y^{2}}}.  \label{2.13}
\end{equation}%
Denote 
\begin{equation}
D^{n+1}=\left\{ \left( x_{1},...,x_{n},y\right) \in \mathbb{R}^{n}\times %
\left[ 0,b\right] \right\} ,\text{ }D_{a}^{n+1}=D^{n+1}\cap \left\{
y>a\right\} .  \label{2.130}
\end{equation}

\textbf{Forward Problem.} \emph{Let (\ref{2.0})-(\ref{2.13}) hold. Find the
function }

$u\left( \mathbf{x},\alpha \right) \in C^{1}\left( \overline{D^{n+1}}\times %
\left[ -d,d\right] \right) $\emph{\ satisfying equation (\ref{2.7}) and the
initial condition}%
\begin{equation}
u(\mathbf{x}_{\alpha },\mathbf{x}_{\alpha })=0\text{ for }\mathbf{x}_{\alpha
}\in \Gamma _{d}.  \label{2.14}
\end{equation}

\textbf{Coefficient Inverse Problem (CIP).} \emph{Let (\ref{2.0})-(\ref{2.13}%
) hold. Let the function }$u\left( \mathbf{x},\alpha \right) \in C^{1}\left( 
\overline{D^{n+1}}\times \left[ -d,d\right] \right) $\emph{\ be the solution
of the Forward Problem. Assume that the coefficient }$a\left( \mathbf{x}%
\right) $\emph{\ of equation (\ref{2.7}) is unknown. Determine the function }%
$a\left( \mathbf{x}\right) ,$\emph{\ assuming that the following function }$%
g\left( \mathbf{x},\alpha \right) $\emph{\ is known:} 
\begin{equation}
g\left( \mathbf{x},\alpha \right) =u\left( \mathbf{x},\alpha \right)
,\forall \mathbf{x}\in \partial \Omega \diagdown \partial _{1}\Omega
,\forall \alpha \in \left( -d,d\right) .  \label{2.15}
\end{equation}

First, we formulate and prove an existence and uniqueness theorem for the
solution of the forward problem. We refer to \cite{SKN,SU} for some other
existence and uniqueness results for the forward problem for equation (\ref%
{2.7}). Unlike Theorem 2.1, the positivity of the function $u$ was not
discussed in these references. On the other hand, we need this property of $u
$ for our numerical method.

\textbf{Theorem 2.1.} \emph{Assume that (\ref{2.0}), (\ref{2.1}), (\ref{2.40}%
)-(\ref{2.50}) and (\ref{2.8})-(\ref{2.13}) hold. Then} \emph{there exists
unique solution }$u\left( \mathbf{x},\alpha \right) \in C^{1}\left( 
\overline{D^{n+1}}\times \left[ -d,d\right] \right) $\emph{\ of equation (%
\ref{2.7}) with the initial condition (\ref{2.14}). In addition, the
following inequality holds:}%
\begin{equation}
u\left( \mathbf{x},\alpha \right) \geq m>0\text{ for }\left( \mathbf{x}%
,\alpha \right) \in \overline{\Omega }\times \left[ -d,d\right] ,
\label{2.16}
\end{equation}%
\begin{equation}
m=\min_{\left( \mathbf{x},\alpha \right) \in \partial _{1}\Omega \times %
\left[ -d.d\right] }\left[ \dint\limits_{L(\mathbf{x},\mathbf{x}_{\alpha
})}f\left( \mathbf{x}\left( s\right) -\mathbf{x}_{\alpha }\right) ds\right]
>0.  \label{2.17}
\end{equation}%
\emph{In addition, there exists a sufficiently large number }$X>A$\emph{\
such that }%
\begin{equation}
u\left( \mathbf{x},\alpha \right) =0\text{ for }\left\vert x_{1}\right\vert
,...,\left\vert x_{n}\right\vert >X,\alpha \in \left( -d,d\right) .
\label{2.18}
\end{equation}

\textbf{Proof}. Equation (\ref{2.7}) can be rewritten as 
\begin{equation}
D_{\nu }u(\mathbf{x},\alpha )+a(\mathbf{x})u(\mathbf{x},\alpha )=\mu _{s}(%
\mathbf{x})\dint\limits_{\Gamma _{d}}K(\mathbf{x},\alpha ,\beta )u(\mathbf{x}%
,\beta )d\beta +f\left( \mathbf{x}-\mathbf{x}_{\alpha }\right) ,
\label{2.19}
\end{equation}%
where $D_{\nu }$ is the directional derivative of $u$ in the direction of
the vector $\nu .$ The first line of (\ref{2.19}) can be treated as the
first order linear ordinary differential operator along the line segment $%
L\left( \mathbf{x},\mathbf{x}_{\alpha }\right) $. Denote%
\begin{equation}
p\left( \mathbf{x},\alpha \right) =\dint\limits_{L(\mathbf{x,x}_{\alpha })}a(%
\mathbf{x}\left( s\right) )ds.  \label{2.189}
\end{equation}%
Obviously, $L(\mathbf{x},\mathbf{x}_{\alpha })=\left\{ \mathbf{z}\left(
t\right) =\left( \alpha +t\left( x_{1}-\alpha \right)
,tx_{2},...,tx_{n},ty\right) ,t\in \left( 0,1\right) \right\} .$ Hence, $%
ds=\left\vert \mathbf{x}-\mathbf{x}_{\alpha }\right\vert dt$ in (\ref{2.189}%
). We obtain 
\begin{equation}
p\left( \mathbf{x},\alpha \right) =\left\vert \mathbf{x}-\mathbf{x}_{\alpha
}\right\vert \dint\limits_{0}^{1}a\left( \alpha +t\left( x_{1}-\alpha
\right) ,tx_{2},...,tx_{n},ty\right) dt.  \label{2.190}
\end{equation}%
Since by (\ref{2.10}) and (\ref{2.12}) $a\in C^{1}\left( \mathbb{R}%
^{n+1}\right) ,$ then (\ref{2.13}), (\ref{2.190}) and elementary
calculations imply that 
\begin{equation}
D_{\nu }p\left( \mathbf{x},\alpha \right) =a\left( \mathbf{x}\right) ,
\label{2.191}
\end{equation}%
where $D_{\nu }$ is the operator of the directional derivative in the
direction of the vector $\nu \left( \mathbf{x},\alpha \right) .$ Consider
now the integration factor $c(\mathbf{x},\alpha ),$ 
\begin{equation}
c(\mathbf{x},\alpha )=e^{p\left( \mathbf{x},\alpha \right) }.  \label{2.20}
\end{equation}%
Then (\ref{2.191}) and (\ref{2.20}) imply 
\begin{equation}
D_{\nu }c(\mathbf{x},\alpha )=a\left( \mathbf{x}\right) c(\mathbf{x},\alpha
).  \label{2.21}
\end{equation}%
Multiply both sides of (\ref{2.19}) by $c(\mathbf{x},\alpha )$. We obtain 
\begin{equation}
c(\mathbf{x,}\alpha )D_{\nu }u(\mathbf{x},\alpha )+c(\mathbf{x},\alpha )a(%
\mathbf{x})u(\mathbf{x},\alpha )=  \label{2.22}
\end{equation}%
\begin{equation*}
=c(\mathbf{x,}\alpha )\mu _{s}(\mathbf{x})\dint\limits_{\Gamma _{d}}K(%
\mathbf{x},\alpha ,\beta )u(\mathbf{x},\beta )d\beta +c(\mathbf{x,}\alpha
)f\left( \mathbf{x}-\mathbf{x}_{\alpha }\right) .
\end{equation*}%
Using (\ref{2.21}), we obtain 
\begin{equation}
cD_{\nu }u+cau=D_{\nu }\left( cu\right) -uD_{\nu }c+cau=D_{\nu }\left(
cu\right) -cau+cau=D_{\nu }\left( cu\right) .  \label{2.23}
\end{equation}%
Since by (\ref{2.5})-(\ref{2.6}), (\ref{2.10}) and (\ref{2.11}) $c(\mathbf{x}%
,\alpha )=1$ for all points $\mathbf{x}$ where $f\left( \mathbf{x}-\mathbf{x}%
_{\alpha }\right) \neq 0,$ then $c(\mathbf{x},\alpha )f\left( \mathbf{x}-%
\mathbf{x}_{\alpha }\right) =f\left( \mathbf{x}-\mathbf{x}_{\alpha }\right)
. $ Hence, (\ref{2.40})-(\ref{2.6}), (\ref{2.11}) and (\ref{2.23}) imply
that problem equation (\ref{2.22}) is equivalent with 
\begin{equation}
D_{\nu }(\left( cu\right) (\mathbf{x,}\alpha ))=c(\mathbf{x,}\alpha )\mu
_{s}(\mathbf{x})\dint\limits_{\Gamma _{d}}K(\mathbf{x},\alpha ,\beta )u(%
\mathbf{x},\beta )d\beta +f\left( \mathbf{x}-\mathbf{x}_{\alpha }\right) .
\label{2.24}
\end{equation}%
Integrating (\ref{2.24}) over the line $L(\mathbf{x},\mathbf{x}_{\alpha })$
and using the initial condition (\ref{2.14}), we obtain \ 
\begin{equation*}
u(\mathbf{x,}\alpha )=\frac{1}{c(\mathbf{x,}\alpha )}\dint\limits_{L(\mathbf{%
x},\mathbf{x}_{\alpha })}c(\mathbf{x}\left( s\right) ,\alpha )\mu _{s}(%
\mathbf{x}\left( s\right) )\left( \dint\limits_{\Gamma _{d}}K(\mathbf{x}%
\left( s\right) ,\alpha ,\beta )u(\mathbf{x}\left( s\right) ,\beta )d\beta
\right) ds+
\end{equation*}%
\begin{equation}
+\frac{1}{c(\mathbf{x,}\alpha )}\dint\limits_{L(\mathbf{x},\mathbf{x}%
_{\alpha })}f\left( \mathbf{x}\left( s\right) -\mathbf{x}_{\alpha }\right)
ds.  \label{2.25}
\end{equation}

Assume now that in (\ref{2.25}) $\mathbf{x}\in \left\{ y\in \left(
0,a\right) \right\} .$ Hence, by (\ref{2.1}) $\mathbf{x\notin }\overline{%
\Omega }.$ Hence, (\ref{2.11}), (\ref{2.189}), (\ref{2.20}) and (\ref{2.25})
imply%
\begin{equation}
u(\mathbf{x,}\alpha )=\dint\limits_{L(\mathbf{x},\mathbf{x}_{\alpha
})}f\left( \mathbf{x}\left( s\right) -\mathbf{x}_{\alpha }\right) ds,\text{ }%
\mathbf{x}\in \left\{ y\in \left( 0,a\right) \right\} .  \label{2.250}
\end{equation}%
This and (\ref{2.2}) imply 
\begin{equation}
u(\mathbf{x,}\alpha )=u_{0}(\mathbf{x,}\alpha )=\dint\limits_{L(\mathbf{x},%
\mathbf{x}_{\alpha })}f\left( \mathbf{x}\left( s\right) -\mathbf{x}_{\alpha
}\right) ds,\text{ }(\mathbf{x,}\alpha )\in \partial _{1}\Omega \times
\left( -d,d\right) ,  \label{2.27}
\end{equation}%
\begin{equation}
u_{0}(\mathbf{x,}\alpha )\geq m,  \label{2.28}
\end{equation}%
where the number $m$ is defined in (\ref{2.17}). In addition, (\ref{2.5}), (%
\ref{2.11}) and (\ref{2.25}) imply that (\ref{2.18}) holds for a
sufficiently large number $X>A.$ It follows from the above construction that
problem (\ref{2.25})- (\ref{2.27}) is equivalent with problem (\ref{2.14}), (%
\ref{2.22}).

Let now $\mathbf{x}\in D_{a}^{n+1},$ where the set $D_{a}^{n+1}$ is defined
in (\ref{2.130}). Hence, similarly with (\ref{2.189}) and (\ref{2.190}), for
any appropriate function $\varphi \left( \mathbf{z}\right) $ 
\begin{equation*}
\dint\limits_{L(\mathbf{x},\mathbf{x}_{\alpha })}\varphi \left( \mathbf{z}%
\left( s\right) \right) ds=\left\vert \mathbf{x}-\mathbf{x}_{\alpha
}\right\vert \dint\limits_{0}^{1}\varphi \left( \alpha +t\left( x_{1}-\alpha
\right) ,tx_{2},...,tx_{n},ty\right) dt=
\end{equation*}%
\begin{equation*}
=\frac{\left\vert \mathbf{x}-\mathbf{x}_{\alpha }\right\vert }{y}%
\dint\limits_{0}^{y}\varphi \left( \alpha +\frac{\left( x_{1}-\alpha \right) 
}{y}z,\frac{x_{2}}{y}z,...,\frac{x_{n}}{y}z,z\right) dz
\end{equation*}%
\begin{equation*}
=\frac{\left\vert \mathbf{x}-\mathbf{x}_{\alpha }\right\vert }{y}%
\dint\limits_{0}^{a}\varphi \left( \alpha +\frac{\left( x_{1}-\alpha \right) 
}{y}z,\frac{x_{2}}{y}z,...,\frac{x_{n}}{y}z,z\right) dz
\end{equation*}%
\begin{equation*}
+\frac{\left\vert \mathbf{x}-\mathbf{x}_{\alpha }\right\vert }{y}%
\dint\limits_{a}^{y}\varphi \left( \alpha +\frac{\left( x_{1}-\alpha \right) 
}{y}z,\frac{x_{2}}{y}z,...,\frac{x_{n}}{y}z,z\right) dz.
\end{equation*}%
Hence, (\ref{2.25}) and (\ref{2.27}) imply%
\begin{equation*}
u(x_{1},x_{2},...,y\mathbf{,}\alpha )=
\end{equation*}%
\begin{equation}
=\frac{\left\vert \mathbf{x}-\mathbf{x}_{\alpha }\right\vert }{yc(\mathbf{x})%
}\dint\limits_{a}^{y}\left( c\mu _{s}\right) (\mathbf{x}\left( z\right)
,\alpha )\left( \dint\limits_{\Gamma _{d}}K(\mathbf{x}\left( z\right)
,\alpha ,\beta )u(\mathbf{x}\left( z\right) ,\beta )d\beta \right) dz+u_{0}(%
\mathbf{x,}\alpha ),  \label{2.29}
\end{equation}%
\begin{equation}
\mathbf{x}\left( z\right) =\left( \alpha +\frac{\left( x_{1}-\alpha \right) 
}{y}z,\frac{x_{2}}{y}z,...,\frac{x_{n}}{y}z,z\right) ,\mathbf{x}\in
D_{a}^{n+1},\alpha \in \left( -d,d\right) .  \label{2.30}
\end{equation}%
Thus, by (\ref{2.18}), we have obtained the $\alpha -$dependent family of
integral equations (\ref{2.18}), (\ref{2.29}), (\ref{2.30}) of the Volterra
type in the bounded domain 
\begin{equation*}
D_{a,b,X}^{n+1}=\left\{ \mathbf{x}=\left( x_{1},...,x_{2},y\right)
:\left\vert x_{1}\right\vert ,\left\vert x_{2}\right\vert ,...\left\vert
x_{n}\right\vert <X,y\in \left( a,b\right) \right\} .
\end{equation*}%
Furthermore, any solution of (\ref{2.29}), (\ref{2.30}) satisfies (\ref{2.18}%
). Since problem (\ref{2.250}), (\ref{2.29}), (\ref{2.30}) is equivalent
with equation (\ref{2.22}) with the initial condition (\ref{2.14}), then it
is sufficient to solve problem (\ref{2.250}), (\ref{2.29}), (\ref{2.30}). It
follows from the well known classical results about Volterra equations that
there exists unique function $u\left( \mathbf{x},\alpha \right) \in C\left( 
\overline{D_{a,b,X}^{n+1}}\times \left[ -d,d\right] \right) $ satisfying (%
\ref{2.29}), (\ref{2.30}), and this function can be obtained via the
following iterative process: 
\begin{equation*}
u_{n}\left( \mathbf{x},\alpha \right) =
\end{equation*}%
\begin{equation}
=\frac{\left\vert \mathbf{x}-\mathbf{x}_{\alpha }\right\vert }{yc(\mathbf{x,}%
\alpha )}\dint\limits_{a}^{y}\left( c\mu _{s}\right) (\mathbf{x}\left(
z\right) ,\alpha )\left( \dint\limits_{\Gamma _{d}}K(\mathbf{x}\left(
z\right) ,\alpha ,\beta )u_{n-1}(\mathbf{x}\left( z\right) ,\beta )d\beta
\right) dz+u_{0}(\mathbf{x,}\alpha ),  \label{2.31}
\end{equation}%
where $n=1,2,...$ It follows from (\ref{2.31}) that for $n=0,1...$ 
\begin{equation}
\left\vert u_{n}\left( \mathbf{x},\alpha \right) \right\vert \leq
\dsum\limits_{k=0}^{n}\frac{\left( M_{1}\left( y-a\right) \right) ^{k}}{k!},%
\text{ }\left( \mathbf{x},\alpha \right) \in \overline{D_{a,b,X}^{n+1}}%
\times \left[ -d,d\right] ,  \label{2.32}
\end{equation}%
where the number $M_{1}=M_{1}\left( a,b,d,X,\left\Vert a\left( x\right)
\right\Vert _{C\left( \overline{\Omega }\right) },\left\Vert K\right\Vert
_{C\left( \overline{\Omega }\times \left[ -d,d\right] ^{2}\right) }\right)
>0 $ depends only on listed parameters. Next, (\ref{2.8}), (\ref{2.27}), (%
\ref{2.28}), (\ref{2.31}) and (\ref{2.32}) imply that estimate (\ref{2.27})
holds. Finally, it follows from (\ref{2.9}), (\ref{2.12}) and (\ref{2.31})
that functions $u_{n}\left( \mathbf{x},\alpha \right) $ can be
differentiated once with respect to $x_{1},...,x_{n},y,\alpha $ and,
similarly with (\ref{2.32}), the following estimates hold for 
\begin{equation}
\left\vert \nabla _{\mathbf{x}}u_{n}\left( \mathbf{x},\alpha \right)
\right\vert ,\left\vert \partial _{\alpha }u_{n}\left( \mathbf{x},\alpha
\right) \right\vert \leq \dsum\limits_{k=0}^{n}\frac{\left( M_{2}\left(
y-a\right) \right) ^{k}}{k!}  \label{2.33}
\end{equation}%
for $\left( \mathbf{x},\alpha \right) \in \overline{D_{a,b,X}^{n+1}}\times %
\left[ -d,d\right] $ and $n=0,1...$ Here the number

$M_{2}=M_{2}\left( a,b,d,X,\left\Vert a\left( x\right) \right\Vert _{C\left( 
\overline{\Omega }\right) },\left\Vert K\right\Vert _{C\left( \overline{%
\Omega }\times \left[ -d,d\right] ^{2}\right) },\left\Vert f\right\Vert
_{C^{1}\left( \left\vert \mathbf{x}\right\vert \leq \varepsilon \right)
}\right) >0$ depends only on listed parameters. Estimates (\ref{2.32}) and (%
\ref{2.33}) together with (\ref{2.18}) imply that we have found the unique
function $u\left( \mathbf{x},\alpha \right) \in C^{1}\left( \overline{D^{n+1}%
}\times \left[ -d,d\right] \right) $ satisfying equation (\ref{2.22}) with
the initial condition (\ref{2.14}). Furthermore, it follows from (\ref{2.8}%
), (\ref{2.28}), (\ref{2.31}) and (\ref{2.32}) that (\ref{2.16}) holds. $\
\square $

\section{Transformation}

\label{sec:3}

\subsection{An integral differential equation without the unknown
coefficient $a\left( \mathbf{x}\right) $}

\label{sec:3.1}

The first step of the convexification is a transformation of the above CIP
to a boundary value problem for a certain PDE of the second order, in which
the unknown coefficient $a\left( \mathbf{x}\right) $ is not involved. By (%
\ref{2.16}) we can introduce a new function $w\left( \mathbf{x},\alpha
\right) ,$ 
\begin{equation}
w\left( \mathbf{x},\alpha \right) =\ln u\left( \mathbf{x},\alpha \right) ,%
\text{ }\left( \mathbf{x},\alpha \right) \in \overline{\Omega }\times \left[
-d,d\right] .  \label{3.1}
\end{equation}%
Substituting (\ref{3.1}) in (\ref{2.7}) and (\ref{2.15}) and using (\ref%
{2.27}), we obtain%
\begin{equation}
\nu (\mathbf{x,}\alpha )\cdot \nabla _{\mathbf{x}}w(\mathbf{x},\alpha
)+a\left( \mathbf{x}\right)  \label{3.2}
\end{equation}%
\begin{equation*}
=e^{-w\left( \mathbf{x},\alpha \right) }\mu _{s}(\mathbf{x}%
)\dint\limits_{\Gamma _{d}}K(\mathbf{x},\alpha ,\beta )e^{w\left( \mathbf{x}%
,\beta \right) }d\beta ,\mathbf{x}\in \Omega ,\alpha \in \left( -d,d\right) ,
\end{equation*}%
\begin{equation}
w(\mathbf{x},\alpha )\mid _{\partial \Omega }=\ln g_{1}\left( \mathbf{x}%
,\alpha \right) ,  \label{3.3}
\end{equation}%
\begin{equation}
g_{1}\left( \mathbf{x},\alpha \right) =\left\{ 
\begin{array}{c}
g\left( \mathbf{x},\alpha \right) ,\mathbf{x}\in \partial \Omega \diagdown
\partial _{1}\Omega ,\alpha \in \left( -d,d\right) , \\ 
u_{0}\left( \mathbf{x},\alpha \right) ,\mathbf{x}\in \partial _{1}\Omega
,\alpha \in \left( -d,d\right) .%
\end{array}%
\right.  \label{3.4}
\end{equation}

Differentiate both sides of (\ref{3.3}) with respect to $\alpha $ and use $%
\partial _{\alpha }a\left( \mathbf{x}\right) \equiv 0.$ We obtain an
integral differential equation with the derivatives up to the second order,%
\begin{equation*}
\nu (\mathbf{x,}\alpha )\cdot \nabla _{\mathbf{x}}w_{\alpha }(\mathbf{x}%
,\alpha )+\partial _{\alpha }\nu (\mathbf{x,}\alpha )\cdot \nabla _{\mathbf{x%
}}w(\mathbf{x},\alpha )
\end{equation*}%
\begin{equation}
=\mu _{s}(\mathbf{x})\frac{\partial }{\partial \alpha }\left[ e^{-w\left( 
\mathbf{x},\alpha \right) }\int_{\Gamma _{d}}K(\mathbf{x},\alpha ,\beta
)e^{w\left( \mathbf{x},\beta \right) }d\beta \right] ,\mathbf{x}\in \Omega
,\alpha \in \left( -d,d\right) .  \label{3.5}
\end{equation}%
We have $\nu (\mathbf{x,}\alpha )=\left( \nu _{1}(\mathbf{x,}\alpha ),\nu
_{2}(\mathbf{x,}\alpha ),...,\nu _{n+1}(\mathbf{x,}\alpha )\right) ,$ where
by (\ref{2.13})%
\begin{equation*}
\nu _{1}(\mathbf{x,}\alpha )=\frac{x_{1}-\alpha }{\sqrt{\left( x_{1}-\alpha
\right) ^{2}+x_{2}^{2}+...+x_{n}^{2}+y^{2}}},
\end{equation*}%
\begin{equation*}
\nu _{k}(\mathbf{x,}\alpha )=\frac{x_{k}}{\sqrt{\left( x_{1}-\alpha \right)
^{2}+x_{2}^{2}+...+x_{n}^{2}+y^{2}}},\text{ }k=2,...,n,
\end{equation*}%
\begin{equation*}
\nu _{n+1}(\mathbf{x,}\alpha )=\frac{y}{\sqrt{\left( x_{1}-\alpha \right)
^{2}+x_{2}^{2}+...+x_{n}^{2}+y^{2}}}.
\end{equation*}%
Hence, (\ref{2.0})-(\ref{2.40}) imply that with a constant $%
M_{3}=M_{3}\left( A,d\right) >0$ depending only on numbers $A,d$ the
following estimates are valid: 
\begin{equation}
\left\vert \frac{\partial _{\alpha }\nu _{k}(\mathbf{x,}\alpha )}{\nu _{n+1}(%
\mathbf{x,}\alpha )}\right\vert ,\left\vert \frac{\nu _{k}(\mathbf{x,}\alpha
)}{\nu _{n+1}(\mathbf{x,}\alpha )}\right\vert \leq \frac{M_{3}}{a},k=1,...,n,%
\mathbf{x\in }\overline{\Omega },\alpha \in \left[ -d,d\right] .  \label{3.6}
\end{equation}

Dividing (\ref{3.5}) by $\nu _{n+1}(\mathbf{x,}\alpha ),$ we obtain%
\begin{equation*}
\partial _{\alpha }w_{y}(\mathbf{x,}\alpha )+\frac{\partial _{\alpha }\nu
_{n+1}(\mathbf{x,}\alpha )}{\nu _{n+1}(\mathbf{x,}\alpha )}w_{y}(\mathbf{x,}%
\alpha )
\end{equation*}%
\begin{equation}
+\dsum\limits_{k=1}^{n}\frac{\nu _{k}(\mathbf{x,}\alpha )}{\nu _{n+1}(%
\mathbf{x,}\alpha )}\partial _{\alpha }w_{x_{k}}(\mathbf{x,}\alpha
)+\dsum\limits_{k=1}^{n}\frac{\partial _{\alpha }\nu _{k}(\mathbf{x,}\alpha )%
}{\nu _{n+1}(\mathbf{x,}\alpha )}w_{x_{k}}(\mathbf{x,}\alpha )  \label{3.8}
\end{equation}%
\begin{equation*}
-\frac{\mu _{s}(\mathbf{x})}{\nu _{n+1}(\mathbf{x,}\alpha )}\frac{\partial }{%
\partial \alpha }\left[ e^{-w\left( \mathbf{x},\alpha \right) }\int_{\Gamma
_{d}}K(\mathbf{x},\alpha ,\beta )e^{w\left( \mathbf{x},\beta \right) }d\beta %
\right] =0,\mathbf{x}\in \Omega ,\alpha \in \left( -d,d\right) ,
\end{equation*}%
\begin{equation}
w(\mathbf{x},\alpha )\mid _{\partial \Omega }=\ln g_{1}\left( \mathbf{x}%
,\alpha \right) .  \label{3.10}
\end{equation}%
The new equation (\ref{3.8}) does not contain the unknown coefficient $%
a\left( \mathbf{x}\right) .$ We need now to solve problem (\ref{3.8}), (\ref%
{3.10}).

\subsection{An orthonormal basis in $L_{2}(-d,d)$}

\label{sec:3.2}

We now describe the orthonormal basis in $L_{2}(-d,d)$ of \cite%
{Klibanov:jiip2017} and \cite[section 6.2.3]{KL}, which was mentioned in
section 1. Consider the set of linearly independent functions $\{\alpha
^{s}e^{\alpha }\}_{s=0}^{\infty }$. This set is complete in the space $%
L_{2}(-d,d)$. Applying the Gram-Schmidt orthonormalization procedure to this
set, we obtain the orthonormal basis $\{\Psi _{s}\left( \alpha \right)
\}_{s=0}^{\infty }$ in $L_{2}(-d,d)$. The function $\Psi _{s}(\alpha )$ has
the form 
\begin{equation}
\Psi _{s}(\alpha )=P_{s}(\alpha )e^{\alpha },\forall s\geq 0,  \label{3.100}
\end{equation}
where $P_{s}(\alpha )$ is a polynomial of the degree $s$.\ Even though the
Gram-Schmidt orthonormalization procedure is unstable for the infinite
number of elements, it is still stable for a non large number of elements,
and this was observed in numerical studies of, e.g. \cite[Chapters 7,10,12]%
{KL}, \cite{KLZ18} and references cited therein. Consider the $N\times N$
matrix $M_{N}=\left( a_{s,k}\right) _{\left( s,k\right) =\left( 0,0\right)
}^{\left( N-1,N-1\right) }$,%
\begin{equation}
a_{s,k}=\mathop{\displaystyle \int}\limits_{-d}^{d}\Psi _{s}^{\prime }\left(
\alpha \right) \Psi _{k}\left( \alpha \right) d\alpha .  \label{3.11}
\end{equation}%
It was proven in \cite{Klibanov:jiip2017}, \cite[section 6.2.3]{KL} that 
\begin{equation}
a_{s,k}=\left\{ 
\begin{array}{c}
1\text{ if }s=k, \\ 
0\text{ if }s<k.%
\end{array}%
\right.  \label{3.13}
\end{equation}%
By (\ref{3.13}) $\det M_{N}=1.$ Hence, the matrix $M_{N}$ is invertible.

\textbf{Remark 3.1}. \emph{The invertibility of the matrix }$M_{N}$\emph{\
is the single most important property of the orthonormal basis }$\{\Psi
_{s}\left( \alpha \right) \}_{s=0}^{\infty }.$\emph{\ Consider, for example
any basis of either classical orthonormal polynomials or orthonormal
trigonometric functions in }$L_{2}(-d,d)$\emph{. Since the first function of
such a basis is a constant, then it follows from (\ref{3.11}) that the first
raw of the analog }$\widetilde{M}_{N}$ \emph{of the matrix }$M_{N}$\emph{\
consists of only zeros. Hence, the matrix }$\widetilde{M}_{N}$\emph{\ is not
invertible.}

\subsection{A coupled system of nonlinear integral differential equations}

\label{sec:3.3}

Represent the function $w(\mathbf{x,}\alpha )$ as truncated Fourier series
with respect to the above basis $\left\{ \Psi _{s}(\alpha )\right\}
_{s=0}^{\infty }$ and also assume that the derivative $w_{\alpha }(\mathbf{x,%
}\alpha )$ can be represented as the sum of the term-by-term derivatives of
that series 
\begin{equation}
w(\mathbf{x,}\alpha )=\dsum\limits_{s=0}^{N-1}w_{s}(\mathbf{x})\Psi
_{s}(\alpha ),\text{ }w_{\alpha }(\mathbf{x,}\alpha
)=\dsum\limits_{n=0}^{N-1}w_{s}(\mathbf{x})\Psi _{n}^{\prime }(\alpha ).
\label{3.19}
\end{equation}%
Therefore, we now need to find the vector function $W(\mathbf{x})$ of
coefficients $w_{s}(\mathbf{x}),$ where 
\begin{equation}
W(\mathbf{x})=\left( w_{0},...,w_{N-1}\right) ^{T}(\mathbf{x}).  \label{3.21}
\end{equation}%
Substituting (\ref{3.19})\emph{\ }in (\ref{3.8}), we obtain 
\begin{equation*}
\dsum\limits_{s=0}^{N-1}\partial _{y}w_{s}(\mathbf{x})\Psi _{s}^{\prime
}(\alpha )+\frac{\partial _{\alpha }\nu _{n+1}(\mathbf{x,}\alpha )}{\nu
_{n+1}(\mathbf{x,}\alpha )}\dsum\limits_{s=0}^{N-1}\partial _{y}w_{sy}(%
\mathbf{x})\Psi _{n}(\alpha )
\end{equation*}%
\begin{equation*}
+\dsum\limits_{s=0}^{N-1}\dsum\limits_{i=1}^{n}\frac{\nu _{i}(\mathbf{x,}%
\alpha )}{\nu _{n+1}(\mathbf{x,}\alpha )}\partial _{x_{i}}w_{s}(\mathbf{x}%
)\Psi _{s}^{\prime }(\alpha )+\dsum\limits_{s=0}^{N-1}\dsum\limits_{i=1}^{n}%
\frac{\partial _{\alpha }\nu _{i}(\mathbf{x,}\alpha )}{\nu _{n+1}(\mathbf{x,}%
\alpha )}\partial _{x_{i}}w_{s}(\mathbf{x})\Psi _{s}(\alpha )
\end{equation*}%
\begin{equation}
-\frac{\mu _{s}(\mathbf{x})}{\nu _{n+1}(\mathbf{x,}\alpha )}\times
\label{3.22}
\end{equation}%
\begin{equation*}
\times \frac{\partial }{\partial \alpha }\left[ \exp \left(
-\dsum\limits_{s=0}^{N-1}w_{s}(\mathbf{x})\Psi _{s}(\alpha )\right)
\dint\limits_{\Gamma _{d}}K(\mathbf{x},\alpha ,\beta )\exp \left(
\dsum\limits_{n=0}^{N-1}w_{s}(\mathbf{x})\Psi _{s}(\beta )\right) d\beta %
\right]
\end{equation*}%
\begin{equation*}
=0,\mathbf{x}\in \Omega ,\alpha \in \left( -d,d\right) .
\end{equation*}%
Using (\ref{3.10}), (\ref{3.19}) and (\ref{3.21}), we obtain the boundary
condition for the vector function $W\left( \mathbf{x}\right) ,$%
\begin{equation}
W\left( \mathbf{x}\right) \mid _{\partial \Omega }=P\left( \mathbf{x}\right)
=\left( p_{0},...,p_{N-1}\right) ^{T}(\mathbf{x}),  \label{3.23}
\end{equation}%
\begin{equation}
p_{s}(\mathbf{x})=\mathop{\displaystyle \int}\limits_{-d}^{d}\ln \left[
g_{1}\left( \mathbf{x},\alpha \right) \right] \Psi _{s}\left( \alpha \right)
d\alpha ,n=0,...,N-1.  \label{3.25}
\end{equation}

Multiply sequentially equation (\ref{3.22}) by functions $\Psi _{k}\left(
\alpha \right) ,k=0,...,N-1$ and integrate with respect to $\alpha \in
\left( -d,d\right) .$ We obtain%
\begin{equation}
\left( M_{N}+A_{n+1}\left( \mathbf{x}\right) \right) W_{y}\left( \mathbf{x}%
\right) +\dsum\limits_{i=1}^{n}A_{i}\left( \mathbf{x}\right) W_{x_{i}}\left( 
\mathbf{x}\right) +F\left( W\left( \mathbf{x}\right) ,\mathbf{x}\right) =0,%
\text{ }\mathbf{x}\in \Omega ,  \label{3.26}
\end{equation}%
where $A_{n+1}\in C_{N^{2}}\left( \overline{\Omega }\right) $ and $A_{i}\in
C_{N^{2}}\left( \overline{\Omega }\right) ,i=1,...,n$ are $N\times N$
matrices, the $N-$D vector function 
\begin{equation}
F\left( s,x\right) \in C^{2}\left( \mathbb{R}^{N+n+1}\right)  \label{3.260}
\end{equation}%
is generated by the operator in the fourth line of (\ref{3.22}). Here, (\ref%
{3.260}) follows from (\ref{2.9}), (\ref{2.10}), (\ref{2.12}), (\ref{3.5})
and (\ref{3.22}). Obviously the vector function $F\left( W\left( \mathbf{x}%
\right) ,\mathbf{x}\right) $ is nonlinear with respect to $W\left( \mathbf{x}%
\right) .$ By (\ref{3.6}) and (\ref{3.7})%
\begin{equation}
\left\Vert A_{i}\left( \mathbf{x}\right) \right\Vert _{C_{N^{2}}\left( 
\overline{\Omega }\right) }\leq \frac{C_{1}}{a},  \label{3.27}
\end{equation}%
where the number $C_{1}=C_{1}\left( A,d,N,b\right) >0$ depends only on
listed parameters.

It follows from (\ref{2.0}) and (\ref{3.27}) that there exists such a number 
$a_{0}>1$ that the matrix 
\begin{equation}
D_{N}\left( \mathbf{x}\right) =\left( M_{N}+A_{n+1}\left( \mathbf{x}\right)
\right) =M_{N}\left( I+M_{N}^{-1}A_{1}\left( \mathbf{x}\right) \right)
\label{3.28}
\end{equation}%
is invertible with the inverse matrix 
\begin{equation}
D_{N}^{-1}\left( \mathbf{x}\right) ,\forall a\geq a_{0}=a_{0}\left(
A,d,N\right) >1,\text{ }\forall \mathbf{x}\in \overline{\Omega },
\label{3.29}
\end{equation}%
where the number $a_{0}\left( A,d,N\right) $ depends only on listed
parameters. Thus, the matrix $D_{N}^{-1}\left( \mathbf{x}\right) $ exists if
the domain $\Omega $ is located sufficiently far from the line of sources $%
\Gamma _{d}.$

\textbf{Remarks 3.2}:

\begin{enumerate}
\item \emph{To work with the convexification method, we need to assume below
that the derivatives }$W_{x_{i}}$ \emph{in (\ref{3.26}) are written in
finite differences, unlike the }$y-$\emph{derivative, and the grid step size 
}$h\geq h_{0}>0$\emph{\ with a fixed constant }$h_{0}$\emph{. The latter
assumption is a practical one since one does not allow the grid step sizes
tend to zero in practical computations.}

\item \emph{The assumptions that the sign \textquotedblleft }$\approx "$%
\emph{\ is replaced with the sign \textquotedblleft }$="$\emph{\ in (\ref%
{3.19}), (\ref{3.22}) and (\ref{3.26}) as well as the assumption of item 1
form our approximate mathematical model, see section 1 for a discussion of
the issue of\ approximate mathematical \ models. Our specific model is
verified computationally in the 2D case in section 6.}
\end{enumerate}

\subsection{Partial finite differences}

\label{sec:3.4}

Let $m>1$ be an integer. Consider $n$ partitions of the interval $\left(
-A,A\right) ,$ see (\ref{2.1}):%
\begin{equation}
-A=x_{i,0}<x_{i,1}<...<x_{i,m}=A,x_{i,j+1}-x_{i,j}=h,j=0,...,m-1,i=1,...,n.
\label{3.30}
\end{equation}%
We assume that 
\begin{equation}
h\geq h_{0}=const.>0.  \label{3.300}
\end{equation}%
Define the discrete subset $\Omega ^{h}\hspace{0.5em}$of the domain $\Omega $
as:%
\begin{equation}
\Omega _{1}^{h}=\left\{ x_{i,j}\right\} _{\left( i,j\right) =\left(
1,0\right) }^{\left( i,j\right) =\left( n,m\right) },\text{ }  \label{3.31}
\end{equation}%
\begin{equation}
\Omega ^{h}=\Omega _{1}^{h}\times \left( a,b\right) =\left\{ \left(
x_{i,j},y\right) :\left\{ x_{i,j}\right\} _{\left( i,j\right) =\left(
1,0\right) }^{\left( i,j\right) =\left( n,m\right) }\in \Omega _{1}^{h},y\in
\left( a,b\right) \right\} .  \label{3.32}
\end{equation}%
We denote $\mathbf{x}^{h}=\left\{ \left( x_{i,j},y\right) :x_{i,j}\in \Omega
_{1}^{h},y\in \left( a,b\right) \right\} .$ By (\ref{2.2})-(\ref{2.4}) and (%
\ref{3.31}), (\ref{3.32}) the boundary $\partial \Omega ^{h}$ of the domain $%
\Omega ^{h}$ is: 
\begin{equation*}
\partial \Omega ^{h}=\partial _{1}\Omega ^{h}\cup \partial _{2}\Omega
^{h}\cup \partial _{3}\Omega ^{h},
\end{equation*}%
\begin{equation*}
\partial _{1}\Omega ^{h}=\Omega _{1}^{h}\times \left\{ y=a\right\} ,\text{ }%
\partial _{2}\Omega ^{h}=\Omega _{1}^{h}\times \left\{ y=b\right\} ,
\end{equation*}%
\begin{equation*}
\partial _{3}\Omega ^{h}=\left\{ \left( x_{i,0},y\right) ,\left(
x_{i,m},y\right) :y\in \left( a,b\right) ,i=1,...,n\right\} .
\end{equation*}%
Let the vector function $Q(\mathbf{x})\in C_{N}^{1}(\overline{\Omega })$.
Denote%
\begin{equation*}
Q_{ij}^{h}\left( y\right) =Q(x_{i,j},y),\hspace{0.5em}i=1,...,n;j=0,...,m;y%
\in \left( a,b\right) ,
\end{equation*}%
\begin{equation*}
Q^{h}(\mathbf{x}^{h})=\left\{ Q^{h}\left(
x_{1,j},x_{2,j},...,x_{n,j},y\right) \right\} ,y\in \left( a,b\right) .
\end{equation*}%
Thus, $Q^{h}(\mathbf{x}^{h})$ is an $N-D$ vector function of discrete
variables $x_{i,j}\in \Omega ^{h}$ and continuous variable $y\in \left(
a,b\right) .$ Note that the boundary terms at $\partial _{3}\Omega ^{h}$ of
this vector function, which correspond to $Q(\mathbf{x})\mid _{\partial
_{3}\Omega ^{h}},$are $\left\{ Q_{i,0}^{h}\left( y\right) \right\} \cup
\left\{ Q_{i,m}^{h}\left( y\right) \right\} ,i=1,...,n.$ For two vector
functions $Q^{\left( 1\right) }(\mathbf{x})=\left( Q_{1}^{\left( 1\right) }(%
\mathbf{x}),...,Q_{N}^{\left( 1\right) }(\mathbf{x})\right) ^{T}$ and $%
Q^{\left( 2\right) }(\mathbf{x})=\left( Q_{1}^{\left( 2\right) }(\mathbf{x}%
),...,Q_{N}^{\left( 2\right) }(\mathbf{x})\right) ^{T}$ their scalar product 
$Q^{\left( 1\right) }(\mathbf{x})\cdot Q^{\left( 2\right) }(\mathbf{x})$ is
defined as the scalar product in $\mathbb{R}^{N}$ and $\left( Q^{\left(
1\right) }(\mathbf{x})\right) ^{2}=Q^{\left( 1\right) }(\mathbf{x})\cdot
Q^{\left( 1\right) }(\mathbf{x}).$ Respectively 
\begin{equation}
Q^{\left( 1\right) h}(\mathbf{x}^{h})\cdot Q^{\left( 2\right) h}(\mathbf{x}%
^{h})=\dsum\limits_{k=1}^{N}\dsum\limits_{\left( i,j\right) =\left(
1,1\right) }^{\left( i,j\right) =\left( 1,m-1\right) }Q_{k}^{\left( 1\right)
h}(x_{i,j},y)Q_{k}^{\left( 2\right) h}(x_{i,j},y),  \label{3.400}
\end{equation}%
\begin{equation}
\left( Q^{h}(\mathbf{x}^{h})\right) ^{2}=Q^{h}(\mathbf{x}^{h})\cdot Q^{h}(%
\mathbf{x}^{h}),\text{ }\left\vert Q^{h}(\mathbf{x}^{h})\right\vert =\sqrt{%
Q^{h}(\mathbf{x}^{h})\cdot Q^{h}(\mathbf{x}^{h})}.  \label{3.401}
\end{equation}%
We will use formulas (\ref{3.400}), (\ref{3.401}) everywhere below without
further mentioning. We exclude $j=0$ and $j=m$ here since we work below with
finite difference derivatives as defined in the next paragraph.

We define finite difference derivatives of $Q^{h}(\mathbf{x}^{h})$ with
respect to $x_{1},...,x_{n}$ only at interior points of the domain $\Omega
^{h}$ with as: 
\begin{equation}
\partial _{x_{k}}Q^{h}(\mathbf{x}^{h})=Q_{x_{k}}^{h}(\mathbf{x}^{h})=\left\{
Q^{h}\left( x_{1j},x_{2j},...,x_{nj},y\right) \right\} _{x_{k}}=
\label{3.40}
\end{equation}%
\begin{equation*}
=\left\{ \frac{Q^{h}\left(
x_{1j},...,x_{k-1,j,}x_{k,j+1},x_{k+1,j}...,x_{n,j},y\right) -Q^{h}\left(
x_{1j},...,x_{k-1,j,}x_{k,j-1},x_{k+1,j}...,x_{n,j},y\right) }{2h}\right\} .
\end{equation*}%
The second line of (\ref{3.40}) should be adjusted in an obvious fashion for 
$k=1$ and $k=n.$ Also, in that line $j=1,...,m-1$. We now define discrete
analogs of spaces $C_{N}\left( \overline{\Omega }\right) $, $H_{N}^{1}\left(
\Omega \right) $ and $L_{N}^{2}\left( \Omega \right) $ as:%
\begin{equation*}
C_{N}^{h}\left( \overline{\Omega ^{h}}\right) =\left\{ Q^{h}(\mathbf{x}%
^{h}):\left\Vert Q^{h}(\mathbf{x}^{h})\right\Vert _{C^{h}\left( \overline{%
\Omega ^{h}}\right) }=\max_{y\in \left[ a,b\right] }\left( \max_{i=\overline{%
1,n},j=\overline{0,m}}\left\vert Q_{ij}^{h}\left( y\right) \right\vert
\right) <\infty \right\} ,
\end{equation*}%
\begin{equation*}
H_{N}^{1,h}\left( \Omega ^{h}\right) =\left\{ 
\begin{array}{c}
Q^{h}(\mathbf{x}^{h}):\left\Vert Q^{h}(\mathbf{x}^{h})\right\Vert
_{H_{N}^{1,h}\left( \Omega ^{h}\right) }^{2}= \\ 
\dsum\limits_{\left( i,j.k\right) =\left( 1,1,1\right) }^{\left(
i,j,k\right) =\left( n,m-1,n\right) }\dint\limits_{a}^{b}\left[ \left(
Q_{ij}^{h}\left( y\right) \right) ^{2}+\left( Q_{x_{k}}^{h}(\mathbf{x}%
^{h})\right) ^{2}+\left( \partial _{y}Q_{ij}^{h}\left( y\right) \right) ^{2}%
\right] dy<\infty%
\end{array}%
\right\} ,
\end{equation*}%
\begin{equation*}
H_{N,0}^{1,h}\left( \Omega ^{h}\right) =\left\{ Q^{h}(\mathbf{x}^{h})\in
H_{N}^{1,h}\left( \Omega ^{h}\right) :Q^{h}(\mathbf{x}^{h})\mid _{\partial
\Omega ^{h}}=0\right\} ,
\end{equation*}%
\begin{equation*}
L_{N}^{2,h}\left( \Omega ^{h}\right) =\left\{ 
\begin{array}{c}
Q^{h}(\mathbf{x}^{h}):\left\Vert Q^{h}(\mathbf{x}^{h})\right\Vert
_{L_{N}^{2,h}\left( \Omega ^{h}\right) }^{2}= \\ 
=\dsum\limits_{\left( i,j\right) =\left( 1,1\right) }^{\left( i,j\right)
=\left( n,m-1\right) }\left\{ \dint\limits_{a}^{b}\left( Q_{ij}^{h}\left(
y\right) \right) ^{2}dy\right\} <\infty%
\end{array}%
\right\} .
\end{equation*}%
By embedding theorem $H_{N}^{1,h}\left( \Omega ^{h}\right) \subset
C_{N}^{h}\left( \overline{\Omega ^{h}}\right) $ and%
\begin{equation}
\left\Vert Q^{h}(\mathbf{x}^{h})\right\Vert _{C_{N}^{h}\left( \overline{%
\Omega ^{h}}\right) }\leq C_{2}\left\Vert Q^{h}(\mathbf{x}^{h})\right\Vert
_{H_{N}^{1,h}\left( \overline{\Omega ^{h}}\right) },\forall Q^{h}\in
H_{N}^{1,h}\left( \Omega ^{h}\right) ,  \label{3.41}
\end{equation}%
where the number $C_{2}=C_{2}\left( h_{0},A,a,b,\Omega \right) >0$ depends
only on listed parameters, and the number $h_{0}$ is defined in (\ref{3.300}%
).

\textbf{Remark 3.3.} \emph{Since we work with finite difference derivatives (%
\ref{3.40}) only at interior points of the discrete domain }$\Omega ^{h}$%
\emph{, then in any differential operator below we use }$W^{h}\left( \mathbf{%
x}^{h}\right) =\left\{ W_{ij}^{h}\left( y\right) \right\} _{\left(
i,j\right) =\left( 1,1\right) }^{\left( i,j\right) =\left( n,m-1\right)
},y\in \left( a,b\right) ,$\emph{\ i.e. boundary points }$x_{i,0}$\emph{\
and }$x_{i,m}$\emph{\ are involved only in finite difference derivatives }$%
W_{i}^{h}\left( \mathbf{x}^{h}\right) $\emph{\ at }$x_{i,1}$\emph{\ and }$%
x_{i,m-1}$\emph{. Points }$x_{i,0}$\emph{\ and }$x_{i,m}$\emph{\ are not
included in (\ref{3.400}) for the same reason.}

Using (\ref{3.30})-(\ref{3.40}), we now rewrite problem (\ref{3.23})-(\ref%
{3.26}) in the form of finite differences with respect to $x_{1},...,x_{n}$
as: 
\begin{equation}
D_{N}^{h}\left( \mathbf{x}^{h}\right) W_{y}^{h}\left( \mathbf{x}^{h}\right)
+\dsum\limits_{i=1}^{n}A_{i}^{h}\left( \mathbf{x}^{h}\right)
W_{x_{i}}^{h}\left( \mathbf{x}^{h}\right) +F^{h}\left( W^{h}\left( \mathbf{x}%
^{h}\right) ,\mathbf{x}^{h}\right) =0,\mathbf{x}^{h}\in \Omega ^{h},
\label{3.42}
\end{equation}%
\begin{equation}
W^{h}\left( \mathbf{x}^{h}\right) \mid _{\partial \Omega ^{h}}=P^{h}\left( 
\mathbf{x}^{h}\right) .  \label{3.43}
\end{equation}%
In (\ref{3.42}) $N\times N$ matrices $A_{i}^{h}\left( \mathbf{x}^{h}\right)
\in C_{N}\left( \overline{\Omega ^{h}}\right) $, the matrix $D_{N}^{h}\left( 
\mathbf{x}^{h}\right) \in C\left( \overline{\Omega ^{h}}\right) $ is the
discrete analog of the matrix $D_{N}\left( \mathbf{x}\right) $ defined in (%
\ref{3.28}) By (\ref{3.27}) and (\ref{3.28}) 
\begin{equation}
\left\Vert D_{N}^{h}\right\Vert _{C_{N^{2}}^{h}\left( \overline{\Omega ^{h}}%
\right) }\leq C_{3},\forall a\geq a_{0}=a_{0}\left( A,d,N\right) >1.
\label{3.430}
\end{equation}%
The vector function $F^{h}\left( W^{h}\left( \mathbf{x}^{h}\right) ,\mathbf{x%
}^{h}\right) $ the discrete analog of the vector function $F$ in (\ref{3.26}%
), and by (\ref{3.260}) $F^{h}\left( W^{h}\left( \mathbf{x}^{h}\right) ,%
\mathbf{x}^{h}\right) $ it is twice continuously differentiable with respect
to components of $W^{h}\left( \mathbf{x}^{h}\right) $. By (\ref{3.29}) there
exists the inverse matrix $\left( D_{N}^{h}\right) ^{-1}\left( \mathbf{x}%
^{h}\right) $ and 
\begin{equation}
\left\Vert \left( D_{N}^{h}\right) ^{-1}\right\Vert _{C_{N^{2}}^{h}\left( 
\overline{\Omega ^{h}}\right) }\leq C_{3},\forall a\geq a_{0}=a_{0}\left(
A,d,N\right) >1.  \label{3.44}
\end{equation}%
Furthermore, (\ref{3.27}) and (\ref{3.44}) imply for $i=1,...,n$ 
\begin{equation}
\left\Vert \left( \left( D_{N}^{h}\right) ^{-1}A_{i}^{h}\right) \left( 
\mathbf{x}^{h}\right) \right\Vert _{C_{N^{2}}^{h}\left( \overline{\Omega ^{h}%
}\right) }\leq C_{3},\forall a\geq a_{0}=a_{0}\left( A,d,N\right) >1,
\label{3.45}
\end{equation}%
Here and everywhere below $C_{3}=C_{3}\left( A,d,N,a,b,h_{0},\left\Vert
K\right\Vert _{C^{1}\left( \overline{\Omega }\times \left[ -d,d\right]
^{2}\right) }\right) >0$ denotes different constants depending only on
listed parameters.

The following formulas are discrete analogs of (\ref{3.19}) and (\ref{3.21}):%
\begin{equation}
w^{h}(\mathbf{x}^{h}\mathbf{,}\alpha )=\dsum\limits_{n=0}^{N-1}w_{n}^{h}(%
\mathbf{x}^{h})\Psi _{n}(\alpha ),\text{ }\partial _{\alpha }w^{h}(\mathbf{x}%
^{h}\mathbf{,}\alpha )=\dsum\limits_{n=0}^{N-1}w_{n}^{h}(\mathbf{x}^{h})\Psi
_{n}^{\prime }(\alpha ),  \label{3.46}
\end{equation}%
\begin{equation}
W^{h}(\mathbf{x}^{h})=\left( w_{0}^{h},...,w_{N-1}^{h}\right) ^{T}(\mathbf{x}%
^{h}).  \label{3.47}
\end{equation}%
Suppose that we have found the vector function $W^{h}(\mathbf{x}^{h})$ in (%
\ref{3.47}). Then, to find the discrete analog $a^{h}\left( \mathbf{x}%
^{h}\right) $ of the unknown coefficient $a\left( \mathbf{x}\right) ,$ we
use (\ref{3.2}) and (\ref{3.46}) as: 
\begin{equation}
a^{h}\left( \mathbf{x}^{h}\right) =-\frac{1}{2d}\dint\limits_{-d}^{d}\nu (%
\mathbf{x}^{h}\mathbf{,}\alpha )\cdot \nabla _{\mathbf{x}^{h}}w^{h}(\mathbf{x%
}^{h},\alpha )d\alpha +  \label{3.48}
\end{equation}%
\begin{equation*}
+\frac{1}{2d}\dint\limits_{-d}^{d}\left( e^{-w^{h}\left( \mathbf{x}%
^{h},\alpha \right) }\mu _{s}(\mathbf{x}^{h})\dint\limits_{-d}^{d}K(\mathbf{x%
}^{h},\alpha ,\beta )e^{w^{h}\left( \mathbf{x}^{h},\beta \right) }d\beta
\right) d\alpha ,\mathbf{x}^{h}\in \Omega ^{h}.
\end{equation*}

\section{Convexification Method for Problem (\protect\ref{3.42}), (\protect
\ref{3.43})}

\label{sec:4}

Let $R>0$ be an arbitrary number. Define the set $B\left( R,P^{h}\right) $ as%
\begin{equation}
B\left( R,P^{h}\right) =\left\{ W^{h}\in H_{N}^{1,h}\left( \Omega
^{h}\right) :W^{h}\left( \mathbf{x}^{h}\right) \mid _{\partial \Omega
^{h}}=P^{h}\left( \mathbf{x}^{h}\right) ,\left\Vert W^{h}\right\Vert
_{H_{N}^{1,h}\left( \Omega ^{h}\right) }<R\right\} ,  \label{4.1}
\end{equation}%
where $P^{h}\left( \mathbf{x}^{h}\right) $ is the boundary condition in (\ref%
{3.43}). Consider matrices $G_{i}^{h}\left( \mathbf{x}^{h}\right) $ and the
vector function $\Phi ^{h}\left( W^{h}\left( \mathbf{x}^{h}\right) ,\mathbf{x%
}^{h}\right) ,$ 
\begin{equation}
G_{i}^{h}\left( \mathbf{x}^{h}\right) =\left( \left( D_{N}^{h}\right)
^{-1}A_{i}^{h}\right) \left( \mathbf{x}^{h}\right) ,\text{ }  \label{4.2}
\end{equation}%
\begin{equation}
\Phi ^{h}\left( W^{h}\left( \mathbf{x}^{h}\right) ,\mathbf{x}^{h}\right)
=\left( D_{N}^{h}\left( \mathbf{x}^{h}\right) \right) ^{-1}F^{h}\left(
W^{h}\left( \mathbf{x}^{h}\right) ,\mathbf{x}^{h}\right) .  \label{4.3}
\end{equation}%
By (\ref{3.260}), (\ref{3.44}), (\ref{3.45}), (\ref{4.1})-(\ref{4.3}) and
the multidimensional analog of Taylor formula \cite{V} 
\begin{equation}
\left\Vert G_{i}^{h}\left( \mathbf{x}^{h}\right) \right\Vert
_{C_{N^{2}}^{h}\left( \overline{\Omega ^{h}}\right) }\leq C_{3},\text{ }%
\forall a\geq a_{0}>1,\text{ }i=1,...,n,  \label{4.4}
\end{equation}%
\begin{equation}
\left\Vert \Phi ^{h}\left( W^{h}\left( \mathbf{x}^{h}\right) ,\mathbf{x}%
^{h}\right) \right\Vert _{C_{N}^{h}\left( \overline{\Omega ^{h}}\right)
}\leq C_{3},\text{ }\forall W^{h}\in \overline{B\left( R,P^{h}\right) },
\label{4.5}
\end{equation}%
\begin{equation}
\Phi ^{h}\left( W_{2}^{h}\left( \mathbf{x}^{h}\right) ,\mathbf{x}^{h}\right)
=\Phi ^{h}\left( W_{1}^{h}\left( \mathbf{x}^{h}\right) ,\mathbf{x}%
^{h}\right) +\Phi _{1}\left( W_{1}^{h}\left( \mathbf{x}^{h}\right) ,\mathbf{x%
}^{h}\right) \left( W_{2}^{h}-W_{1}^{h}\right) \left( \mathbf{x}^{h}\right)
\label{4.6}
\end{equation}%
\begin{equation*}
+\Phi _{2}^{h}\left( W_{1}^{h}\left( \mathbf{x}^{h}\right) ,W_{2}^{h}\left( 
\mathbf{x}^{h}\right) ,\mathbf{x}^{h}\right) ,\text{ }\forall
W_{1}^{h},W_{2}^{h}\in \overline{B\left( R,P^{h}\right) },\forall \mathbf{x}%
^{h}\in \overline{\Omega ^{h}},
\end{equation*}%
where the vector function $\Phi _{1}\left( W_{1}^{h}\left( \mathbf{x}%
^{h}\right) ,\mathbf{x}^{h}\right) $ is independent on $W_{2}^{h}\left( 
\mathbf{x}^{h}\right) ,$ the vector function $\Phi _{2}^{h}\left(
W_{1}^{h}\left( \mathbf{x}^{h}\right) ,W_{2}^{h}\left( \mathbf{x}^{h}\right)
,\mathbf{x}^{h}\right) $ is nonlinear with respect to $\left(
W_{2}^{h}-W_{1}^{h}\right) \left( \mathbf{x}^{h}\right) ,$ both these vector
functions are continuous with respect to their variables for $\mathbf{x}%
^{h}\in \Omega ^{h}$ and the following estimates hold for all $%
W_{1}^{h},W_{2}^{h}\in \overline{B\left( R,P^{h}\right) }$ and for all $%
\mathbf{x}^{h}\in \Omega ^{h}$%
\begin{equation}
\left\vert \Phi _{1}\left( W_{1}^{h}\left( \mathbf{x}^{h}\right) ,\mathbf{x}%
^{h}\right) \right\vert \leq C_{3},\text{ }\forall W_{1}^{h}\in \overline{%
B\left( R,P^{h}\right) },\forall \mathbf{x}^{h}\in \Omega ^{h},  \label{4.60}
\end{equation}%
\begin{equation}
\left\vert \Phi _{2}^{h}\left( W_{1}^{h}\left( \mathbf{x}^{h}\right)
,W_{2}^{h}\left( \mathbf{x}^{h}\right) ,\mathbf{x}^{h}\right) \right\vert
\leq C_{3}\left( W_{2}^{h}-W_{1}^{h}\right) ^{2}\left( \mathbf{x}^{h}\right)
,\text{ }\forall \mathbf{x}^{h}\in \Omega ^{h},  \label{4.61}
\end{equation}%
also, see (\ref{3.401}).

\textbf{Lemma 4.1}. \emph{Let }$A$\emph{\ be a }$k\times k$\emph{\ matrix
which has the inverse }$A^{-1}.$\emph{\ Then there exists a number }$\mu
=\mu \left( A\right) >0$\emph{\ such that }$\left\Vert Ax\right\Vert
^{2}\geq \mu \left\Vert x\right\Vert ^{2},\forall x\in \mathbb{R}^{k},$ 
\emph{where }$\left\Vert \cdot \right\Vert $\emph{\ is} \emph{the euclidean
norm.}

We omit the proof since this lemma is well known.

\textbf{Corollary 4.1}. \emph{The following inequality holds:}%
\begin{equation*}
\left( D_{N}^{h}\left( \mathbf{x}^{h}\right) W_{y}^{h}\left( \mathbf{x}%
^{h}\right) +\dsum\limits_{i=1}^{n}A_{i}^{h}\left( \mathbf{x}^{h}\right)
W_{x_{i}}^{h}\left( \mathbf{x}^{h}\right) \right) ^{2}\geq
\end{equation*}%
\begin{equation*}
\geq C_{3}\left( W_{y}^{h}\left( \mathbf{x}^{h}\right)
+\dsum\limits_{i=1}^{n}G_{i}^{h}\left( \mathbf{x}^{h}\right)
W_{x_{i}}^{h}\left( \mathbf{x}^{h}\right) \right) ^{2},\forall \mathbf{x}%
^{h}\in \Omega ^{h}.
\end{equation*}

\textbf{Proof}. Denote 
\begin{equation}
Y\left( \mathbf{x}^{h}\right) =W_{y}^{h}\left( \mathbf{x}^{h}\right)
+\dsum\limits_{i=1}^{n}G_{i}^{h}\left( \mathbf{x}^{h}\right)
W_{x_{i}}^{h}\left( \mathbf{x}^{h}\right) +\Phi ^{h}\left( W^{h}\left( 
\mathbf{x}^{h}\right) ,\mathbf{x}^{h}\right) .  \label{4.62}
\end{equation}%
We have 
\begin{equation}
\left( D_{N}^{h}\left( \mathbf{x}^{h}\right) W_{y}^{h}\left( \mathbf{x}%
^{h}\right) +\dsum\limits_{i=1}^{n}A_{i}^{h}\left( \mathbf{x}^{h}\right)
W_{x_{i}}^{h}\left( \mathbf{x}^{h}\right) \right) ^{2}=\left(
D_{N}^{h}\left( \mathbf{x}^{h}\right) Y\left( \mathbf{x}^{h}\right) \right)
^{2}.  \label{4.63}
\end{equation}%
The rest of the proof follows immediately from (\ref{4.2}), (\ref{4.3}) and
Lemma 4.1. $\square $

Introduce the following weighted cost functional $J_{\lambda }\left(
W^{h}\right) $:%
\begin{equation}
J_{\lambda }\left( W^{h}\right) =\left\Vert \left(
D_{N}^{h}W_{y}^{h}+\dsum\limits_{i=1}^{n}A_{i}^{h}W_{x_{i}}^{h}+F^{h}\left(
W^{h}\left( \mathbf{x}^{h}\right) ,\mathbf{x}^{h}\right) \right) e^{\lambda
y}\right\Vert _{L_{N}^{2,h}\left( \Omega ^{h}\right) }^{2}.  \label{4.7}
\end{equation}

\textbf{Minimization\ Problem. }\emph{Minimize functional (\ref{4.7}) on the
set} $\overline{B\left( R,P^{h}\right) }.$

Theorems 4.1-4.5 are our analytical results about the functional $J_{\lambda
,\gamma }\left( W^{h}\right) .$

\textbf{Theorem 4.1 }(Carleman estimate). \emph{Assume that the number }$%
a\geq a_{0},$\emph{\ as in (\ref{3.44}). Then there exists a sufficiently
large number }$\lambda _{0}=\lambda _{0}\left( A,d,N,a,b,h_{0}\right) \geq 1$%
\emph{\ depending only on listed parameters such that the following Carleman
estimate holds }%
\begin{equation}
\left\Vert \left( W_{y}^{h}\left( \mathbf{x}^{h}\right)
+\dsum\limits_{i=1}^{n}G_{i}^{h}\left( \mathbf{x}^{h}\right)
W_{x_{i}}^{h}\left( \mathbf{x}^{h}\right) \right) e^{\lambda y}\right\Vert
_{L_{N}^{2,h}\left( \Omega ^{h}\right) }^{2}\geq  \label{4.8}
\end{equation}%
\begin{equation*}
\geq \left\Vert W_{y}^{h}e^{\lambda y}\right\Vert _{L_{N}^{2,h}\left( \Omega
^{h}\right) }^{2}+\lambda ^{2}\left\Vert W^{h}e^{\lambda y}\right\Vert
_{L_{N}^{2,h}\left( \Omega ^{h}\right) }^{2},\text{ }\forall W^{h}\in
H_{N,0}^{1,h}\left( \Omega ^{h}\right) ,\forall \lambda \geq \lambda _{0}.
\end{equation*}%
\emph{\ }

\textbf{Theorem 4.2} (the central analytical result). \emph{The following
three assertions hold:}

\bigskip 1. \emph{The functional }$J_{\lambda }\left( W^{h}\right) $\emph{\
in (\ref{4.7}) has the Fr\'{e}chet derivative }$J_{\lambda }^{\prime }\left(
W^{h}\right) \in H_{N,0}^{1,h}\left( \Omega ^{h}\right) $\emph{\ at any
point }$W^{h}\in \overline{B\left( R,P^{h}\right) }$\emph{\ and for any
value of the parameter }$\lambda \geq 0.$\emph{\ This function satisfies the
Lipschitz condition}%
\begin{equation}
\left\Vert J_{\lambda }^{\prime }\left( W_{1}^{h}\right) -J_{\lambda
}^{\prime }\left( W_{2}^{h}\right) \right\Vert _{H_{N}^{1,h}\left( \Omega
^{h}\right) }\leq C_{4}\left\Vert W_{1}^{h}-W_{2}^{h}\right\Vert
_{H_{N}^{1,h}\left( \Omega ^{h}\right) },\forall W_{1}^{h},W_{2}^{h}\in 
\overline{B\left( R,P^{h}\right) }  \label{4.80}
\end{equation}%
\emph{for all }$\lambda \geq 0,$ \emph{where the number }$C_{4}>0$\emph{\
depends on the same parameters as ones in }$C_{3}$\emph{\ as well as on }$%
\lambda .$

\emph{Assume that the number }$a\geq a_{0},$\emph{\ as in (\ref{3.44}). Then:%
}

\emph{2. There exists a sufficiently large number }$\lambda _{1}$\emph{\ }%
\begin{equation}
\lambda _{1}=\lambda _{1}\left( R,A,d,N,a,b,h_{0}\right) \geq \lambda
_{0}\geq 1  \label{4.81}
\end{equation}%
\emph{depending only on listed parameters such that the functional }$%
J_{\lambda }\left( W^{h}\right) $\emph{\ in (\ref{4.7}) is strictly convex
on the set }$\overline{B\left( R,P^{h}\right) },$\emph{\ i.e. the following
inequality holds:}%
\begin{equation}
J_{\lambda }\left( W_{2}^{h}\right) -J_{\lambda }\left( W_{1}^{h}\right)
-J_{\lambda }^{\prime }\left( W_{1}^{h}\right) \left(
W_{2}^{h}-W_{1}^{h}\right) \geq C_{3}e^{2\lambda a}\left\Vert
W_{2}^{h}-W_{1}^{h}\right\Vert _{H_{N}^{1,h}\left( \Omega ^{h}\right) }^{2},
\label{4.9}
\end{equation}%
\begin{equation}
\forall \lambda \geq \lambda _{1},\forall W_{1}^{h},W_{2}^{h}\in \overline{%
B\left( R,P^{h}\right) }.  \label{4.90}
\end{equation}

\emph{3. For each }$\lambda \geq \lambda _{1}$\emph{\ there exists unique
minimizer }$W_{\min ,\lambda }^{h}\in \overline{B\left( R,P^{h}\right) }$%
\emph{\ of the functional }$J_{\lambda }\left( W^{h}\right) $\emph{\ on the
set }$\overline{B\left( R,P^{h}\right) }$\emph{\ . Furthermore, the
following inequality holds: }%
\begin{equation}
J_{\lambda }^{\prime }\left( W_{\min ,\lambda }^{h}\right) \left(
W^{h}-W_{\min ,\lambda }^{h}\right) \geq 0,\text{ }\forall W^{h}\in 
\overline{B\left( R,P^{h}\right) }.  \label{4.10}
\end{equation}

Theorem 4.3 follows immediately from (\ref{4.9}) and (\ref{4.90}). This is a
certain uniqueness result for our CIP, which is obtained as a by-product. A
further discussion of the uniqueness issue is outside of the scope of this
paper.

\textbf{Theorem 4.3}. \emph{Assume that the number }$a\geq a_{0},$\emph{\ as
in (\ref{3.44}). Then there exists at most one pair of functions }$\left(
W^{h},a^{h}\right) \in H_{N}^{1,h}\left( \Omega ^{h}\right) \times
L_{N}^{2,h}\left( \Omega ^{h}\right) $ \emph{satisfying conditions (\ref%
{3.42}), (\ref{3.48}).}

We now estimate the accuracy of the minimizer $W_{\min ,\lambda }^{h}$
depending on the level of the noise $\delta >0$ in the data. Following the
concept of Tikhonov for ill-posed problems \cite{T}, we assume the existence
of the exact solution 
\begin{equation}
W^{h\ast }\in B\left( R,P^{h\ast }\right)  \label{4.100}
\end{equation}%
of problem (\ref{3.42})-(\ref{3.43}) with the exact, i.e. noiseless data $%
P^{h\ast },$ i.e. for $\mathbf{x}^{h}\in \Omega ^{h}$ 
\begin{equation}
D_{N}^{h}\left( \mathbf{x}^{h}\right) W_{y}^{h\ast }\left( \mathbf{x}%
^{h}\right) +\dsum\limits_{i=1}^{n}A_{i}^{h}\left( \mathbf{x}^{h}\right)
W_{x_{i}}^{h\ast }\left( \mathbf{x}^{h}\right) +F^{h}\left( W^{h\ast }\left( 
\mathbf{x}^{h}\right) ,\mathbf{x}^{h}\right) =0,  \label{4.11}
\end{equation}%
\begin{equation}
W^{h\ast }\left( \mathbf{x}^{h}\right) \mid _{\partial \Omega ^{h}}=P^{h\ast
}\left( \mathbf{x}^{h}\right) .  \label{4.13}
\end{equation}%
Suppose that there exists a vector function $S^{h}\in H_{N}^{1,h}\left(
\Omega ^{h}\right) $ such that 
\begin{equation}
S^{h}\left( \mathbf{x}^{h}\right) \mid _{\partial \Omega ^{h}}=P^{h}\left( 
\mathbf{x}^{h}\right) ,\text{ }\left\Vert S^{h}\right\Vert
_{H_{N}^{1,h}\left( \Omega ^{h}\right) }\,<R.  \label{4.14}
\end{equation}%
Let $S^{h\ast }\in H_{N}^{1,h}\left( \Omega ^{h}\right) $ be such a vector
function that 
\begin{equation}
S^{h\ast }\left( \mathbf{x}^{h}\right) \mid _{\partial \Omega ^{h}}=P^{h\ast
}\left( \mathbf{x}^{h}\right) ,\text{ }\left\Vert S^{h\ast }\right\Vert
_{H_{N}^{1,h}\left( \Omega ^{h}\right) }\,<R.  \label{4.15}
\end{equation}%
We assume that 
\begin{equation}
\left\Vert S^{h}-S^{h\ast }\right\Vert _{H_{N}^{1,h}\left( \Omega
^{h}\right) }<\delta .  \label{4.16}
\end{equation}

\textbf{Theorem 4.4}. \emph{Assume that the number }$a\geq a_{0},$\emph{\ as
in (\ref{3.44}). Suppose that conditions (\ref{4.11})-(\ref{4.16}) hold.
Also, consider the number }$\lambda _{2},$\emph{\ }%
\begin{equation}
\lambda _{2}=\lambda _{1}\left( 2R,A,d,N,a,b,h_{0}\right) ,  \label{4.160}
\end{equation}%
\emph{where }$\lambda _{1}\left( R,A,d,N,a,b,h_{0}\right) $\emph{\ is the
number in (\ref{4.81}). Let }$W_{\min ,\lambda _{2}}^{h}$\emph{\ be the
minimizer of functional (\ref{4.7}) on the set }$\overline{B\left(
R,P^{h}\right) },$\emph{\ which was found in Theorem 4.2. Let }$\alpha \in
\left( 0,R\right) $ \emph{be a number. Suppose that (\ref{4.100}) is
replaced with}%
\begin{equation}
W^{h\ast }\in B\left( R-\alpha ,P^{h\ast }\right) \text{ and }C_{3}\delta
<\alpha .  \label{4.161}
\end{equation}%
\emph{Then the vector function }$W_{\min ,\lambda _{2}}^{h}$\emph{\ belongs
to the open \ set }$B\left( R,P^{h}\right) $\emph{\ and the following
accuracy estimate holds:}%
\begin{equation}
\left\Vert W_{\min ,\lambda _{2}}^{h}-W^{h\ast }\right\Vert
_{H_{N}^{1,h}\left( \Omega ^{h}\right) }\leq C_{3}\delta .  \label{4.17}
\end{equation}%
\emph{\ }

Consider now the gradient descent method of the minimization of functional 
\emph{(\ref{4.7}) on the set }$\overline{B\left( R,P^{h}\right) }.$ Let $%
W_{0}^{h}\in B\left( R/3,P^{h}\right) $ be an arbitrary point of this set.
We treat it as the starting point of the latter method. The sequence of this
method is:%
\begin{equation}
W_{n}^{h}=W_{n-1}^{h}-\gamma J_{\lambda _{2}}^{\prime }\left(
W_{n-1}^{h}\right) ,n=1,2,...,  \label{4.18}
\end{equation}%
where $\gamma >0$ is a small number and $\lambda _{2}$ is the same as in (%
\ref{4.160}). Note that since by Theorem 4.2 functions $J_{\lambda
_{2}}^{\prime }\left( W_{n-1}^{h}\right) \in H_{N,0}^{1,h}\left( \Omega
^{h}\right) $, then all vector functions $W_{n}^{h}$ \ have the same
boundary conditions $P^{h}.$

\textbf{Theorem 4.5.} \emph{Let conditions of Theorem 4.4 hold, except that (%
\ref{4.161}) is replaced with} 
\begin{equation}
W^{h\ast }\in B\left( \left( R-\alpha \right) /3,P^{h\ast }\right) \text{
and }C_{3}\delta <\alpha /3.  \label{4.19}
\end{equation}%
\emph{Then there exists a sufficiently small number }$\gamma >0$\emph{\ and
a number }$\theta =\theta \left( \gamma \right) \in \left( 0,1\right) $\emph{%
\ such that in (\ref{4.18}) all functions }$W_{n}^{h}\in B\left(
R,P^{h}\right) $\emph{\ and the following convergence estimates hold:}%
\begin{equation}
\left\Vert W_{n}^{h}-W_{\min ,\lambda _{2}}^{h}\right\Vert
_{H_{N}^{1,h}\left( \Omega ^{h}\right) }\leq \theta ^{n}\left\Vert
W_{0}^{h}-W_{\min ,\lambda _{2}}^{h}\right\Vert _{H_{N}^{1,h}\left( \Omega
^{h}\right) },  \label{4.20}
\end{equation}%
\begin{equation}
\left\Vert W_{n}^{h}-W^{h\ast }\right\Vert _{H_{N}^{1,h}\left( \Omega
^{h}\right) }\leq C_{3}\delta +\theta ^{n}\left\Vert W_{0}^{h}-W_{\min
,\lambda _{2}}^{h}\right\Vert _{H_{N}^{1,h}\left( \Omega ^{h}\right) },
\label{4.21}
\end{equation}%
\begin{equation}
\left\Vert a_{n}^{h}-a^{h\ast }\right\Vert _{L_{N}^{2,h}\left( \Omega
^{h}\right) }\leq C_{3}\delta +\theta ^{n}\left\Vert W_{0}^{h}-W_{\min
,\lambda _{2}}^{h}\right\Vert _{H_{N}^{1,h}\left( \Omega ^{h}\right) },
\label{4.22}
\end{equation}%
\emph{where }$a_{n}^{h}\left( \mathbf{x}^{h}\right) $\emph{\ and }$%
a_{n}^{h\ast }\left( \mathbf{x}^{h}\right) $\emph{\ are functions which are
obtained from }$W_{n}^{h}$\emph{\ and }$W^{h\ast }$\emph{\ respectively via (%
\ref{3.46}) and (\ref{3.48}).}

\textbf{Remarks 4.1}:

\begin{enumerate}
\item \emph{Estimates (\ref{4.20})-(\ref{4.22}) guarantee the global
convergence of the gradient descent method (\ref{4.18}) since }$R>0$\emph{\
is an arbitrary number and the starting point }$W_{0}^{h}$\emph{\ is an
arbitrary point of the set }$B\left( R/3,P^{h}\right) $\emph{, see section 1
for our definition of the global convergence. Note that any gradient-like
method converges only locally for a non-convex functional, i.e. it needs a
good first guess about the correct solution. }

\item \emph{Although the above results are valid only for sufficiently large
values of the parameter }$\lambda ,$\emph{\ we have established
computationally in section 6 that }$\lambda =5$\emph{\ works quite well. It
was computationally established in a number of works on the convexification
for a variety of CIPs that values }$\lambda \in \left[ 1,3\right] $\emph{\
work well numerically, see, e.g. \cite[Chapters 7-10]{KL}, \cite{KLZ18} and
references cited therein.}

\item \emph{An analogy to the issue of item 2 is that any asymptotic theory
basically states that if a parameter }$X$\emph{\ is sufficiently
large/small, then a certain `nice' formula provides a good approximation for
something. In practice, however, given a specific problem with specific
values of its parameters, only numerical experiments can clarify which
exactly values of }$X$\emph{\ are so sufficiently large/small that this
`nice' formula indeed provides that good approximation. }
\end{enumerate}

\section{Proofs}

\label{sec:5}

We use in this section (\ref{3.400})-(\ref{3.40}) and Remark 3.3 without
further mentioning.

\subsection{Proof of Theorem 4.1}

\label{sec:5.1}

In this proof, $W^{h}\in H_{N,0}^{1,h}\left( \Omega ^{h}\right) $ is an
arbitrary function. It follows from (\ref{3.300}) and (\ref{3.40}) that 
\begin{equation}
\left\Vert W_{x_{i}}^{h}e^{\lambda y}\right\Vert _{L_{N}^{2,h}\left( \Omega
^{h}\right) }^{2}\leq C_{3}\left\Vert W^{h}e^{\lambda y}\right\Vert
_{L_{N}^{2,h}\left( \Omega ^{h}\right) }^{2},\text{ }\forall W^{h}\in
H_{N,0}^{1,h}\left( \Omega ^{h}\right) .  \label{5.100}
\end{equation}%
By (\ref{4.4}) and (\ref{5.100})%
\begin{equation}
\left\Vert \left( \dsum\limits_{i=1}^{n}G_{i}^{h}\left( \mathbf{x}%
^{h}\right) W_{x_{i}}^{h}\left( \mathbf{x}^{h}\right) \right) e^{\lambda
y}\right\Vert _{L_{N}^{2,h}\left( \Omega ^{h}\right) }^{2}\leq
C_{3}\left\Vert W^{h}e^{\lambda y}\right\Vert .  \label{5.2}
\end{equation}%
Using Cauchy-Schwarz inequality and (\ref{5.2}), we obtain%
\begin{equation*}
\left\Vert \left( W_{y}^{h}\left( \mathbf{x}^{h}\right)
+\dsum\limits_{i=1}^{n}G_{i}^{h}\left( \mathbf{x}^{h}\right)
W_{x_{i}}^{h}\left( \mathbf{x}^{h}\right) \right) e^{\lambda y}\right\Vert
_{L_{N}^{2,h}\left( \Omega ^{h}\right) }^{2}\geq
\end{equation*}%
\begin{equation}
\geq C_{3,1}\left\Vert W_{y}^{h}e^{\lambda y}\right\Vert _{L_{N}^{2,h}\left(
\Omega ^{h}\right) }^{2}-C_{3,2}\left\Vert W^{h}e^{\lambda y}\right\Vert
_{L_{N}^{2,h}\left( \Omega ^{h}\right) }^{2}.  \label{5.3}
\end{equation}%
It makes sense to use two numbers $C_{3,1},C_{3,2}>0$ here, both depend on
the same parameters as ones in $C_{3}.$ Consider now the first term in the
second line of (\ref{5.3}). Introduce a new vector function $V^{h}\left( 
\mathbf{x}^{h}\right) =W^{h}\left( \mathbf{x}^{h}\right) e^{\lambda y}.$
Then 
\begin{equation*}
W^{h}\left( \mathbf{x}^{h}\right) =V^{h}\left( \mathbf{x}^{h}\right)
e^{-\lambda y},W_{y}^{h}\left( \mathbf{x}^{h}\right) =\left( V_{y}^{h}\left( 
\mathbf{x}^{h}\right) -\lambda V^{h}\left( \mathbf{x}^{h}\right) \right)
e^{-\lambda y},
\end{equation*}%
\begin{equation*}
\left( W_{y}^{h}\left( \mathbf{x}^{h}\right) \right) ^{2}e^{2\lambda
y}=\left( V_{y}^{h}\left( \mathbf{x}^{h}\right) -\lambda V^{h}\left( \mathbf{%
x}^{h}\right) \right) ^{2}\geq -2\lambda V_{y}^{h}\left( \mathbf{x}%
^{h}\right) V^{h}\left( \mathbf{x}^{h}\right) +\lambda ^{2}\left(
V^{h}\left( \mathbf{x}^{h}\right) \right) ^{2}
\end{equation*}%
\begin{equation*}
=\partial _{y}\left( -\lambda V^{h}\left( \mathbf{x}^{h}\right) \right)
^{2}+\lambda ^{2}\left( V^{h}\left( \mathbf{x}^{h}\right) \right) ^{2}.
\end{equation*}%
Hence, 
\begin{equation*}
\left\Vert W_{y}^{h}e^{\lambda y}\right\Vert _{L_{N}^{2,h}\left( \Omega
^{h}\right) }^{2}\geq \dsum\limits_{\left( i,j\right) =\left( 1,1\right)
}^{\left( i,j\right) =\left( n,m-1\right) }\left[ -\lambda \left(
W^{h}\left( \mathbf{x}^{h}\right) \right) ^{2}\mid _{y=b}+\lambda \left(
W^{h}\left( \mathbf{x}^{h}\right) \right) ^{2}\mid _{y=a}\right] +
\end{equation*}%
\begin{equation}
+\lambda ^{2}\left\Vert W^{h}e^{\lambda y}\right\Vert _{L_{N}^{2,h}\left(
\Omega ^{h}\right) }^{2}.  \label{5.4}
\end{equation}%
Since $W^{h}\left( \mathbf{x}^{h}\right) \mid _{y=b}=0,$ $\forall W^{h}\in
H_{N,0}^{1,h}\left( \Omega ^{h}\right) ,$ then (\ref{5.4}) implies%
\begin{equation}
\left\Vert W_{y}^{h}e^{\lambda y}\right\Vert _{L_{N}^{2,h}\left( \Omega
^{h}\right) }^{2}\geq \lambda ^{2}\left\Vert W^{h}e^{\lambda y}\right\Vert
_{L_{N}^{2,h}\left( \Omega ^{h}\right) }^{2}.  \label{5.5}
\end{equation}%
Adding the term $\left\Vert W_{y}^{h}e^{\lambda y}\right\Vert
_{L_{N}^{2,h}\left( \Omega ^{h}\right) }^{2}$ to both sides of (\ref{5.5})
and then dividing by \textquotedblleft 2", we obtain for all $\lambda >0$ 
\begin{equation}
\left\Vert W_{y}^{h}e^{\lambda y}\right\Vert _{L_{N}^{2,h}\left( \Omega
^{h}\right) }^{2}\geq \frac{1}{2}\left\Vert W_{y}^{h}e^{\lambda
y}\right\Vert _{L_{N}^{2,h}\left( \Omega ^{h}\right) }^{2}+\frac{\lambda ^{2}%
}{2}\left\Vert W^{h}e^{\lambda y}\right\Vert _{L_{N}^{2,h}\left( \Omega
^{h}\right) }^{2}.  \label{5.50}
\end{equation}%
Hence, taking $\lambda _{0}^{2}\geq 4C_{3,2}/C_{3,1},\lambda \geq \lambda
_{0}$ and using (\ref{5.3}) and (\ref{5.50}), we obtain the target estimate (%
\ref{4.8}) of this theorem. $\square $

\subsection{Proof of Theorem 4.2}

\label{sec:5.2}

Let $W_{1}^{h}\left( \mathbf{x}^{h}\right) ,W_{2}^{h}\left( \mathbf{x}%
^{h}\right) \in \overline{B\left( R,P^{h}\right) }$ be two arbitrary points.
Denote $v^{h}\left( \mathbf{x}^{h}\right) =W_{2}^{h}\left( \mathbf{x}%
^{h}\right) -W_{1}^{h}\left( \mathbf{x}^{h}\right) .$ By (\ref{4.1})%
\begin{equation}
v^{h}\in H_{N,0}^{1,h}\left( \Omega ^{h}\right) .  \label{5.6}
\end{equation}%
By (\ref{4.7}) 
\begin{equation}
J_{\lambda }\left( W_{2}^{h}\right) =J_{\lambda }\left(
W_{1}^{h}+v^{h}\right) =\dsum\limits_{\left( i,j\right) =\left( 1,1\right)
}^{\left( i,j\right) =\left( n,m-1\right) }\times  \label{5.7}
\end{equation}%
\begin{equation*}
\times \dint\limits_{a}^{b}\left( D_{N}^{h}\left(
W_{1y}^{h}+v_{y}^{h}\right) +\dsum\limits_{i=1}^{n}A_{i}^{h}\left(
W_{1x_{i}}^{h}+v_{x_{i}}\right) +F^{h}\left( \left( W^{h}+v^{h}\right) ,%
\mathbf{x}^{h}\right) \right) ^{2}e^{2\lambda y}dy.
\end{equation*}%
Consider the integrand in the second line of (\ref{5.7}) without, however,
the term $e^{2\lambda y}.$ By (\ref{4.2}), (\ref{4.3}), (\ref{4.62}) and (%
\ref{4.63}) 
\begin{equation}
\left( D_{N}^{h}\left( W_{1y}^{h}+v_{y}^{h}\right)
+\dsum\limits_{i=1}^{n}A_{i}^{h}\left( \mathbf{x}^{h}\right) \left(
W_{1x_{i}}^{h}+v_{x_{i}}^{h}\right) +F^{h}\left( W^{h}+v^{h},\mathbf{x}%
^{h}\right) \right) ^{2}  \label{5.8}
\end{equation}%
\begin{equation*}
=\left[ D_{N}^{h}\left( \left( W_{1y}^{h}+v_{y}^{h}\right)
+\dsum\limits_{i=1}^{n}G_{i}^{h}\left( W_{1x_{i}}^{h}+v_{x_{i}}^{h}\right)
+\Phi ^{h}\left( W_{1}^{h}+v^{h},\mathbf{x}^{h}\right) \right) \right] ^{2}.
\end{equation*}%
Consider the term $\Phi ^{h}\left( W_{1}^{h}+v^{h},\mathbf{x}^{h}\right) .$
By (\ref{4.6})-(\ref{4.61}) we have for all $\mathbf{x}^{h}\in \Omega ^{h}$ 
\begin{equation}
\Phi ^{h}\left( W_{1}^{h}+v^{h},\mathbf{x}^{h}\right) =\Phi ^{h}\left(
W_{1}^{h},\mathbf{x}^{h}\right) +\Phi _{1}\left( W_{1}^{h},\mathbf{x}%
^{h}\right) v^{h}\left( \mathbf{x}^{h}\right) +\Phi _{2}^{h}\left(
W_{1}^{h},W_{1}^{h}\mathbf{+}v^{h},\mathbf{x}^{h}\right) ,  \label{5.9}
\end{equation}%
\begin{equation}
\left\vert \Phi _{1}\left( W_{1}^{h},\mathbf{x}^{h}\right) \right\vert \leq
C_{3},\text{ }\left\vert \Phi _{2}^{h}\left( W_{1}^{h},W_{1}^{h}\mathbf{+}%
v^{h},\mathbf{x}^{h}\right) \right\vert \leq C_{3}\left( v^{h}\left( \mathbf{%
x}^{h}\right) \right) ^{2}\text{.}  \label{5.10}
\end{equation}%
It follows from (\ref{5.9}) and (\ref{5.10}) that the second line of (\ref%
{5.8}) is: 
\begin{equation*}
\left[ D_{N}^{h}\left( \left( W_{1y}^{h}+v_{y}^{h}\right)
+\dsum\limits_{i=1}^{n}G_{i}^{h}\left( W_{1x_{i}}^{h}+v_{x_{i}}^{h}\right)
+\Phi ^{h}\left( W_{1}^{h}+v^{h},\mathbf{x}^{h}\right) \right) \right] ^{2}
\end{equation*}%
\begin{equation*}
=\left[ D_{N}^{h}\left(
W_{1y}^{h}+\dsum\limits_{i=1}^{n}G_{i}^{h}W_{1x_{i}}^{h}+\Phi ^{h}\left(
W_{1}^{h},\mathbf{x}^{h}\right) \right) \right] ^{2}
\end{equation*}%
\begin{equation*}
+2D_{N}^{h}\left(
W_{1y}^{h}+\dsum\limits_{i=1}^{n}G_{i}^{h}W_{1x_{i}}^{h}+\Phi ^{h}\left(
W_{1}^{h},\mathbf{x}^{h}\right) \right) \cdot
\end{equation*}%
\begin{equation}
\cdot D_{N}^{h}\left(
v_{y}^{h}+\dsum\limits_{i=1}^{n}G_{i}^{h}v_{x_{i}}^{h}+\Phi _{1}\left(
W_{1}^{h},\mathbf{x}^{h}\right) v^{h}\left( \mathbf{x}^{h}\right) \right)
\label{5.11}
\end{equation}%
\begin{equation*}
+2D_{N}^{h}\left(
W_{1y}^{h}+\dsum\limits_{i=1}^{n}G_{i}^{h}W_{1x_{i}}^{h}+\Phi ^{h}\left(
W_{1}^{h},\mathbf{x}^{h}\right) \right) \cdot \left( D_{N}^{h}\Phi
_{2}^{h}\left( W_{1}^{h},W_{1}^{h}\mathbf{+}v^{h},\mathbf{x}^{h}\right)
\right)
\end{equation*}%
\begin{equation*}
+\left[ D_{N}^{h}\left(
v_{y}^{h}+\dsum\limits_{i=1}^{n}G_{i}^{h}v_{x_{i}}^{h}+\Phi _{1}\left(
W_{1}^{h},\mathbf{x}^{h}\right) v^{h}\left( \mathbf{x}^{h}\right) +\Phi
_{2}^{h}\left( W_{1}^{h},W_{1}^{h}\mathbf{+}v^{h},\mathbf{x}^{h}\right)
\right) \right] ^{2}.
\end{equation*}%
The linear with respect to $v^{h}$ term in (\ref{5.11}) is the scalar
product of the third and fourth lines of this equality. We denote the latter
as $Lin\left( v^{h}\left( \mathbf{x}^{h}\right) ,\mathbf{x}^{h}\right) .$
The sum of fifth and sixth lines of (\ref{5.11}) contains only nonlinear
terms with respect to $v^{h}.$ We denote this sum as $Nonlin\left(
v^{h}\left( \mathbf{x}^{h}\right) ,\mathbf{x}^{h}\right) .$ It follows from
Corollary 4.1, (\ref{3.430}), (\ref{4.2})-(\ref{4.61}), (\ref{5.10}) and (%
\ref{5.11}) and Cauchy-Schwarz inequality that%
\begin{equation}
\left\vert Nonlin\left( v^{h}\left( \mathbf{x}^{h}\right) ,\mathbf{x}%
^{h}\right) \right\vert \geq C_{3,3}\left(
v_{y}^{h}+\dsum\limits_{i=1}^{n}G_{i}^{h}v_{x_{i}}^{h}\right)
^{2}-C_{3,4}\left( v^{h}\left( \mathbf{x}^{h}\right) \right) ^{2},\forall 
\mathbf{x}^{h}\in \overline{\Omega ^{h}}.  \label{5.12}
\end{equation}%
Just as in the proof of Theorem 4.1, it again makes sense to use two
positive constants $C_{3,3},C_{3,4}>0$ here, both depend on the same
parameters as ones involved in $C_{3}.$ Taking into account (\ref{5.7}), (%
\ref{5.8}) and the second line of (\ref{5.11}), we obtain%
\begin{equation}
J_{\lambda }\left( W_{2}^{h}\right) =J_{\lambda }\left(
W_{1}^{h}+v^{h}\right) =J_{\lambda }\left( W_{1}^{h}\right) +  \label{5.13}
\end{equation}%
\begin{equation*}
+\dsum\limits_{\left( i,j\right) =\left( 1,1\right) }^{\left( i,j\right)
=\left( n,m-1\right) }\dint\limits_{a}^{b}Lin\left( v^{h}\left( \mathbf{x}%
^{h}\right) ,\mathbf{x}^{h}\right) e^{2\lambda y}dy+\dsum\limits_{\left(
i,j\right) =\left( 1,1\right) }^{\left( i,j\right) =\left( n,m-1\right)
}\dint\limits_{a}^{b}Nonlin\left( v^{h}\left( \mathbf{x}^{h}\right) ,\mathbf{%
x}^{h}\right) e^{2\lambda y}dy.
\end{equation*}%
It follows from (\ref{4.4})-(\ref{4.6}), (\ref{5.6}) and from the third and
fourth lines of (\ref{5.11}), which form $Lin\left( v^{h}\left( \mathbf{x}%
^{h}\right) ,\mathbf{x}^{h}\right) ,$ that the first term in the second line
of (\ref{5.13}) is a bounded linear functional with respect to $v^{h}$
mapping $H_{N,0}^{1,h}\left( \Omega ^{h}\right) $ in $\mathbb{R}.$ We denote
this functional as $\widetilde{J}_{\lambda }\left( W_{1}^{h}\right) \left(
v^{h}\right) .$ By Riesz theorem there exists unique function $J_{\lambda
}^{\prime }\left( W_{1}^{h}\right) \in H_{N,0}^{1,h}\left( \Omega
^{h}\right) $ such that 
\begin{equation}
\widetilde{J}_{\lambda }\left( W_{1}^{h}\right) \left( v^{h}\right) =\left(
J_{\lambda }^{\prime }\left( W_{1}^{h}\right) ,v^{h}\right) ,\forall
v^{h}\in H_{N,0}^{1,h}\left( \Omega ^{h}\right) ,  \label{5.14}
\end{equation}%
where $\left( ,\right) $ is the scalar product in $H_{N,0}^{1,h}\left(
\Omega ^{h}\right) .$ It follows from (\ref{3.44}), (\ref{3.45}), (\ref{4.4}%
)-(\ref{4.60}) and the fifth and sixth lines of (\ref{5.11}) combined with (%
\ref{5.10}) that for $\left\Vert v^{h}\right\Vert _{H_{N}^{1,h}\left( \Omega
^{h}\right) }<1$ 
\begin{equation*}
\left\vert \dsum\limits_{\left( i,j\right) =\left( 1,1\right) }^{\left(
i,j\right) =\left( n,m-1\right) }\dint\limits_{a}^{b}Nonlin\left(
v^{h}\left( \mathbf{x}^{h}\right) ,\mathbf{x}^{h}\right) e^{2\lambda
y}dy\right\vert \leq C_{3}e^{2\lambda b}\left\Vert v^{h}\right\Vert
_{H_{N}^{1,h}\left( \Omega ^{h}\right) }^{2}.\text{ }
\end{equation*}%
Hence, (\ref{5.13}) implies that 
\begin{equation*}
\frac{J_{\lambda }\left( W_{1}^{h}+v^{h}\right) -J_{\lambda }\left(
W_{1}^{h}\right) -\left( J_{\lambda }^{\prime }\left( W_{1}^{h}\right)
,v^{h}\right) }{\left\Vert v^{h}\right\Vert _{H_{N}^{1,h}\left( \Omega
^{h}\right) }}=O\left( \left\Vert v^{h}\right\Vert _{H_{N}^{1,h}\left(
\Omega ^{h}\right) }\right) ,\text{ as }\left\Vert v^{h}\right\Vert
_{H_{N}^{1,h}\left( \Omega ^{h}\right) }\rightarrow 0.
\end{equation*}%
Thus, $J_{\lambda }^{\prime }\left( W_{1}^{h}\right) :H_{N,0}^{1,h}\left(
\Omega ^{h}\right) \rightarrow \mathbb{R}$ is the Fr\'{e}chet derivative of
the functional $J_{\lambda }\left( W^{h}\right) $ at the point $W_{1}^{h}.$
We omit the proof of estimate (\ref{4.80}) since this proof is completely
similar with the proof of Theorem 3.1 of \cite{Bak}.

We now prove (\ref{4.9}), (\ref{4.90}). Using (\ref{5.12})-(\ref{5.14}), we
obtain%
\begin{equation*}
J_{\lambda }\left( W_{1}^{h}+v^{h}\right) -J_{\lambda }\left(
W_{1}^{h}\right) -\left( J_{\lambda }^{\prime }\left( W_{1}^{h}\right)
,v^{h}\right) \geq
\end{equation*}%
\begin{equation}
\geq C_{3,3}\left\Vert \left(
v_{y}^{h}+\dsum\limits_{i=1}^{n}G_{i}^{h}v_{x_{i}}^{h}\right) e^{\lambda
y}\right\Vert _{L_{N}^{2,h}\left( \Omega ^{h}\right) }^{2}-C_{3,4}\left\Vert
ve^{\lambda y}\right\Vert _{L_{N}^{2,h}\left( \Omega ^{h}\right) }^{2}.
\label{5.15}
\end{equation}%
Let the number $\lambda _{0}$ be the one chosen in Theorem 4.1. Recalling (%
\ref{5.6}) and using Theorem 4.1, we estimate the second line of (\ref{5.15}%
) for all $\lambda \geq \lambda _{0},$%
\begin{equation}
C_{3,3}\left\Vert \left(
v_{y}^{h}+\dsum\limits_{i=1}^{n}G_{i}^{h}v_{x_{i}}^{h}\right) e^{\lambda
y}\right\Vert _{L_{N}^{2,h}\left( \Omega ^{h}\right) }^{2}-C_{3,4}\left\Vert
v^{h}e^{\lambda y}\right\Vert _{L_{N}^{2,h}\left( \Omega ^{h}\right)
}^{2}\geq  \label{5.149}
\end{equation}%
\begin{equation*}
\geq C_{3}\left( \left\Vert v_{y}^{h}e^{\lambda y}\right\Vert
_{L_{N}^{2,h}\left( \Omega ^{h}\right) }^{2}+\lambda ^{2}\left\Vert
v^{h}e^{\lambda y}\right\Vert _{L_{N}^{2,h}\left( \Omega ^{h}\right)
}^{2}\right) -C_{3,4}\left\Vert v^{h}e^{\lambda y}\right\Vert
_{L_{N}^{2,h}\left( \Omega ^{h}\right) }^{2}.
\end{equation*}%
Choosing $\lambda _{1}\geq \lambda _{0}$ so large that $\lambda
_{1}^{2}C_{3}/2>2C_{3,4}$ and using (\ref{5.149}), we obtain%
\begin{equation}
C_{3,3}\left\Vert \left(
v_{y}^{h}+\dsum\limits_{i=1}^{n}G_{i}^{h}v_{x_{i}}^{h}\right) e^{\lambda
y}\right\Vert _{L_{N}^{2,h}\left( \Omega ^{h}\right) }^{2}-C_{3,4}\left\Vert
v^{h}e^{\lambda y}\right\Vert _{L_{N}^{2,h}\left( \Omega ^{h}\right)
}^{2}\geq  \label{5.150}
\end{equation}%
\begin{equation*}
\geq C_{3}\left( \left\Vert v_{y}e^{\lambda y}\right\Vert
_{L_{N}^{2,h}\left( \Omega ^{h}\right) }^{2}+\lambda ^{2}\left\Vert
ve^{\lambda y}\right\Vert _{L_{N}^{2,h}\left( \Omega ^{h}\right)
}^{2}\right) .
\end{equation*}%
Combining (\ref{5.15}) with (\ref{5.150}), we obtain (\ref{4.9}), (\ref{4.90}%
).

Finally, given (\ref{4.9}) and (\ref{4.90}), existence and uniqueness of the
minimizer $W_{\min ,\lambda }^{h}$ of the functional $J_{\lambda }\left(
W^{h}\right) $ on the set $\overline{B\left( R,P^{h}\right) }$ as well as
estimate (\ref{4.10}) follow immediately from a combination of Lemma 2.1
with Theorem \ 2.1 of \cite{Bak}. $\ \ \square $

\subsection{Proof of Theorem 4.4}

\label{sec:5.3}

Denote 
\begin{equation}
B_{0}\left( 2R\right) =\left\{ V^{h}\in H_{N,0}^{1,h}\left( \Omega
^{h}\right) :\left\Vert V^{h}\right\Vert _{H_{N}^{1,h}\left( \Omega
^{h}\right) }<2R\right\} .  \label{5.16}
\end{equation}%
Consider the vector functions $V^{h},V^{h\ast }$%
\begin{equation}
V^{h}=W^{h}-S^{h},\forall W^{h}\in B\left( R,P^{h}\right) ;\text{ }V^{h\ast
}=W^{h\ast }-S^{h\ast }.  \label{5.17}
\end{equation}%
By (\ref{4.100}), (\ref{4.14}), (\ref{4.15}), (\ref{5.16}) and (\ref{5.17}) 
\begin{equation}
V^{h},V^{h\ast }\in B_{0}\left( 2R\right) .  \label{5.18}
\end{equation}%
Consider now the functional $I_{\lambda }\left( V^{h}\right) ,$%
\begin{equation}
I_{\lambda }\left( V^{h}\right) :B_{0}\left( 2R\right) \rightarrow \mathbb{R}%
,I_{\lambda }\left( V^{h}\right) =J_{\lambda }\left( V^{h}+S^{h}\right) .
\label{5.19}
\end{equation}%
An obvious analog of Theorem 4.2 holds for $I_{\lambda }\left( V^{h}\right)
. $ However, it follows from (\ref{5.18}) that we need to replace $R$ with $%
2R$ in (\ref{4.81}) in this case, i.e. we need to use now the number $%
\lambda _{2}$ in (\ref{4.160}). Let $V_{\min ,\lambda _{2}}^{h}$ be the
minimizer of $I_{\lambda _{2}}\left( V^{h}\right) $ on the set $B_{0}\left(
2R\right) .$

Consider $I_{\lambda }\left( V^{h\ast }\right) =J_{\lambda }\left( V^{h\ast
}+S^{h}\right) .$ By (\ref{4.9}), (\ref{4.90}) and (\ref{5.19}) 
\begin{equation*}
I_{\lambda _{2}}\left( V^{h\ast }\right) -I_{\lambda _{2}}\left( V_{\min
,\lambda _{2}}^{h}\right) -I_{\lambda _{2}}^{\prime }\left( V_{\min ,\lambda
_{2}}^{h}\right) \left( V^{h\ast }-V_{\min ,\lambda _{2}}^{h}\right) \geq
C_{3}\left\Vert V^{h\ast }-V_{\min ,\lambda _{2}}^{h}\right\Vert
_{H_{N}^{1,h}\left( \Omega ^{h}\right) }^{2}.
\end{equation*}%
By (\ref{4.10}) $-I_{\lambda _{2}}^{\prime }\left( V_{\min ,\lambda
_{2}}^{h}\right) \left( V^{h\ast }-V_{\min ,\lambda _{2}}^{h}\right) \leq 0.$
Since $-I_{\lambda _{2}}\left( V_{\min ,\lambda _{2}}^{h}\right) \leq 0$ as
well, then the latter estimate implies 
\begin{equation}
\left\Vert V^{h\ast }-V_{\min ,\lambda _{2}}^{h}\right\Vert
_{H_{N}^{1,h}\left( \Omega ^{h}\right) }^{2}\leq C_{3}I_{\lambda _{2}}\left(
V^{h\ast }\right) .  \label{5.21}
\end{equation}%
Next, by (\ref{5.17}) and (\ref{5.19})%
\begin{equation}
I_{\lambda _{2}}\left( V^{h\ast }\right) =J_{\lambda _{2}}\left( V^{h\ast
}+S^{h}\right) =J_{\lambda _{2}}\left( W^{h\ast }+\left( S^{h}-S^{h\ast
}\right) \right) .  \label{5.22}
\end{equation}%
By (\ref{4.11}) the right hand side of (\ref{4.7}) equals zero if $W^{h}$ is
replaced with $W^{h\ast }.$ Hence, using (\ref{4.2})-(\ref{4.61}), (\ref%
{4.16}) and (\ref{5.17}), we obtain $J_{\lambda _{2}}\left( W^{h\ast
}+\left( S^{h}-S^{h\ast }\right) \right) \leq C_{3}\delta ^{2}.$ Hence, (\ref%
{5.21}) and (\ref{5.22}) lead to:%
\begin{equation}
\left\Vert V^{h\ast }-V_{\min ,\lambda _{2}}^{h}\right\Vert
_{H_{N}^{1,h}\left( \Omega ^{h}\right) }\leq C_{3}\delta .  \label{5.23}
\end{equation}%
Using again (\ref{4.16}), (\ref{5.17}), triangle inequality and (\ref{5.23}%
), we obtain 
\begin{equation}
\left\Vert W^{h\ast }-\widetilde{W}_{\min ,\lambda _{2}}^{h}\right\Vert
_{H_{N}^{1,h}\left( \Omega ^{h}\right) }\leq C_{3}\delta .  \label{5.24}
\end{equation}%
Here $\widetilde{W}_{\min ,\lambda _{2}}^{h}=V_{\min ,\lambda
_{2}}^{h}-S^{h}.$ By (\ref{4.161}), (\ref{5.24}) and the triangle inequality 
$\left\Vert \widetilde{W}_{\min ,\lambda _{2}}^{h}\right\Vert
_{H_{N}^{1,h}\left( \Omega ^{h}\right) }<R.$ Hence, (\ref{4.1}), (\ref{4.14}%
), (\ref{5.16}) and (\ref{5.18}) imply that 
\begin{equation}
\widetilde{W}_{\min ,\lambda _{2}}^{h}\in B\left( R,P^{h}\right) .
\label{5.25}
\end{equation}

Now, let $W_{\min ,\lambda _{2}}^{h}$ be the minimizer of the functional $%
J_{\lambda _{2}}\left( W^{h}\right) $ on the set $\overline{B\left(
R,P^{h}\right) },$ which is claimed by Theorem 4.2. Let $Q_{\min ,\lambda
_{2}}^{h}=W_{\min ,\lambda _{2}}^{h}-S^{h}.$ Then $Q_{\min ,\lambda
_{2}}^{h}\in B_{0}\left( 2R\right) .$ Hence, 
\begin{equation*}
J_{\lambda _{2}}\left( \widetilde{W}_{\min ,\lambda _{2}}^{h}\right)
=J_{\lambda _{2}}\left( V_{\min ,\lambda _{2}}^{h}+G\right) \leq J_{\lambda
_{2}}\left( Q_{\min ,\lambda _{2}}^{h}+G\right) =\min_{\overline{B\left(
R,P^{h}\right) }}J_{\lambda _{2}}\left( W^{h}\right) .
\end{equation*}%
Hence, by (\ref{5.25}) $\widetilde{W}_{\min ,\lambda _{2}}^{h}$ is also a
minimizer of the functional $J_{\lambda _{2}}\left( W^{h}\right) $ on the
set $\overline{B\left( R,P^{h}\right) }.$ However, since such a minimizer is
unique by Theorem 4.2, then $\widetilde{W}_{\min ,\lambda _{2}}^{h}=W_{\min
,\lambda _{2}}^{h}\in B\left( R,P^{h}\right) .$ Hence, estimate (\ref{5.24})
holds when $\widetilde{W}_{\min ,\lambda _{2}}^{h}$ is replaced with $%
W_{\min ,\lambda _{2}}^{h}.$ The latter immediately implies (\ref{4.17}).$\
\square $

\subsection{Proof of Theorem 4.5}

\label{sec:5.4}

Consider again the minimizer $W_{\min ,\lambda _{2}}^{h}$ of the functional $%
J_{\lambda _{2}}\left( W^{h}\right) $ on the set $\overline{B\left(
R,P^{h}\right) }.$ Then (\ref{4.19}) and Theorem 4.4 imply that $W_{\min
,\lambda _{2}}^{h}\in B\left( R/3,P^{h}\right) .$ To prove (\ref{4.20}), we
now modify the proof of Theorem 6 of \cite{SAR}: to adapt it for our
specific case.

We set $\lambda =\lambda _{2},$ where $\lambda =\lambda _{2}$ is defined in (%
\ref{4.160}). It obviously follows from (\ref{4.9}) that%
\begin{equation}
\left( J_{\lambda _{2}}^{\prime }\left( W_{1}^{h}\right) -J_{\lambda
_{2}}^{\prime }\left( W_{2}^{h}\right) \right) \left(
W_{1}^{h}-W_{2}^{h}\right) \geq C_{3}\left\Vert
W_{2}^{h}-W_{1}^{h}\right\Vert _{H_{N}^{1,h}\left( \Omega ^{h}\right) }^{2}
\label{5.26}
\end{equation}%
for all $W_{1}^{h},W_{2}^{h}\in B\left( R,P^{h}\right) .$ Consider the
operator $T\left( W^{h}\right) :B\left( R,P^{h}\right) \rightarrow
H_{N}^{1,h}\left( \Omega ^{h}\right) ,$ 
\begin{equation}
T\left( W^{h}\right) =W^{h}-\gamma J_{\lambda _{2}}^{\prime }\left(
W^{h}\right) .  \label{5.27}
\end{equation}%
Using (\ref{4.80}) and (\ref{5.26}), a contraction property of the operator $%
X\left( W^{h}\right) ,$ it was proven in \cite[Theorem 2.1]{Bak} and in \cite%
[page 97]{KL} that there exists a sufficiently small number $\gamma >0$\ and
a number $\theta =\theta \left( \gamma \right) \in \left( 0,1\right) $\ such
that%
\begin{equation}
\left\Vert T\left( W_{1}^{h}\right) -T\left( W_{2}^{h}\right) \right\Vert
_{H_{N}^{1,h}\left( \Omega ^{h}\right) }\leq \theta \left\Vert
W_{1}^{h}-W_{2}^{h}\right\Vert _{H_{N}^{1,h}\left( \Omega ^{h}\right) }
\label{5.28}
\end{equation}%
for all $W_{1}^{h},W_{2}^{h}\in B\left( R,P^{h}\right) .$ However, it was
not proven in \cite{Bak,KL} that the operator $T$ maps $B\left(
R,P^{h}\right) $ in itself.

Since $W_{\min ,\lambda _{2}}^{h}$ is an interior point of the set $B\left(
R,P^{h}\right) ,$ then $J_{\lambda _{2}}^{\prime }\left( W_{\min ,\lambda
_{2}}^{h}\right) =0.$ Hence, by (\ref{5.27}) 
\begin{equation}
W_{\min ,\lambda _{2}}^{h}=T\left( W_{\min ,\lambda _{2}}^{h}\right) .
\label{5.29}
\end{equation}
Consider now the vector function $W_{1}^{h}=T\left( W_{0}^{h}\right)
=W_{0}^{h}-\gamma J_{\lambda _{2}}^{\prime }\left( W_{0}^{h}\right) .$ Since
by Theorem 4.2 the function $J_{\lambda _{2}}^{\prime }\left(
W_{0}^{h}\right) \in H_{N,0}^{1,h}\left( \Omega ^{h}\right) ,$ then $%
W_{1}^{h}\mid _{\partial \Omega ^{h}}=W_{0}^{h}\mid _{\partial \Omega
^{h}}=P^{h}.$ Next by (\ref{5.28}), (\ref{5.29}) and triangle inequality 
\begin{equation*}
\left\Vert W_{1}^{h}-W_{\min ,\lambda _{2}}^{h}\right\Vert
_{H_{N}^{1,h}\left( \Omega ^{h}\right) }=\left\Vert T\left( W_{0}^{h}\right)
-T\left( W_{\min ,\lambda _{2}}^{h}\right) \right\Vert _{H_{N}^{1,h}\left(
\Omega ^{h}\right) }\leq
\end{equation*}%
\begin{equation}
\leq \theta \left\Vert W_{0}^{h}-W_{\min ,\lambda _{2}}^{h}\right\Vert
_{H_{N}^{1,h}\left( \Omega ^{h}\right) }\leq \frac{2R}{3}\theta <\frac{2R}{3}%
.  \label{5.30}
\end{equation}%
Hence, (\ref{4.161}), (\ref{4.17}) and triangle inequality imply $\left\Vert
W_{1}^{h}\right\Vert _{H_{N}^{1,h}\left( \Omega ^{h}\right) }<R.$ Hence, by (%
\ref{4.1}) $W_{1}^{h}\in B\left( R,P^{h}\right) .$ Similarly by (\ref{5.30}) 
\begin{equation*}
\left\Vert W_{2}^{h}-W_{\min ,\lambda _{2}}^{h}\right\Vert
_{H_{N}^{1,h}\left( \Omega ^{h}\right) }\leq \theta \left\Vert
W_{1}^{h}-W_{\min ,\lambda _{2}}^{h}\right\Vert _{H_{N}^{1,h}\left( \Omega
^{h}\right) }\leq
\end{equation*}%
\begin{equation*}
\leq \theta ^{2}\left\Vert W_{0}^{h}-W_{\min ,\lambda _{2}}^{h}\right\Vert
_{H_{N}^{1,h}\left( \Omega ^{h}\right) }<\frac{2R}{3}.
\end{equation*}%
Hence, applying again triangle inequality, we obtain $W_{2}^{h}\in B\left(
R,P^{h}\right) .$ Continuing this process, we obtain (\ref{4.20}).

To prove (\ref{4.21}), we note that by triangle inequality, (\ref{4.17}) and
(\ref{4.20}) 
\begin{equation*}
\left\Vert W_{n}^{h}-W^{h\ast }\right\Vert _{H_{N}^{1,h}\left( \Omega
^{h}\right) }\leq \left\Vert W_{\min ,\lambda _{2}}^{h}-W^{h\ast
}\right\Vert _{H_{N}^{1,h}\left( \Omega ^{h}\right) }+\left\Vert
W_{n}^{h}-W_{\min ,\lambda _{2}}^{h}\right\Vert _{H_{N}^{1,h}\left( \Omega
^{h}\right) }\leq
\end{equation*}%
\begin{equation*}
\leq C_{3}\delta +\theta ^{n}\left\Vert W_{0}^{h}-W_{\min ,\lambda
_{2}}^{h}\right\Vert _{H_{N}^{1,h}\left( \Omega ^{h}\right) }.
\end{equation*}%
The proof of estimate (\ref{4.22}) follows immediately from (\ref{3.46})-(%
\ref{3.48}) and (\ref{4.21}). \ $\square $

\section{Numerical Studies}

\label{sec:6}

In order not to introduce new and complicated notations, we slightly abuse
below some notations of the previous sections. Nevertheless, the substance
is always clear from the context presented below.

\subsection{Numerical implementation}

\label{sec:6.1}

We have conducted our numerical studies in the 2-D case. Below $\mathbf{x}%
=\left( x,y\right) $ and, according to (\ref{2.1}) and (\ref{2.40}) 
\begin{equation}
\Omega =\{\mathbf{x}:x\in \left( -A,A\right) ,y\in \left( a,b\right)
,A=1/2,a=1,b=2,  \label{6.0}
\end{equation}%
\begin{equation*}
\Gamma _{d}=\left\{ \mathbf{x}_{\alpha }=(\alpha ,0):\alpha \in \lbrack
-d,d]\right\} ,d=1/2.
\end{equation*}%
As to the kernel $K(\mathbf{x},\alpha ,\beta )$ of the integral operator in (%
\ref{2.7}), we work below with the 2D Henyey-Greenstein function \cite{Heino}%
: 
\begin{equation}
K(\mathbf{x},\alpha ,\beta )=H(\alpha ,\beta )=\frac{1}{2d}\left[ \frac{%
1-g^{2}}{1+g^{2}-2g\cos (\alpha -\beta )}\right] ,\quad g=\frac{1}{2}.
\label{6.2}
\end{equation}%
Here, $g=1/2$ means an anisotropic scattering, which is half ballistic with $%
g=0$ an half isotropic scattering with $g=1$ \cite{HT1,HT2,HT3}. We take the
same function $f\left( \mathbf{x}\right) $ as one in (\ref{2.5}), (\ref{2.50}%
) with $\varepsilon =0.05.$

We assume that 
\begin{equation}
\mu _{s}(\mathbf{x})=\left\{ 
\begin{array}{ll}
5, & \quad \mathbf{x}\in \Omega , \\ 
0, & \quad \mathbf{x}\in \mathbb{R}^{2}\setminus \Omega .%
\end{array}%
\right.  \label{6.3}
\end{equation}%
Given (\ref{6.3}), we use formula (\ref{2.10}) for the coefficient function $%
a(\mathbf{x})$ and we take in this formula 
\begin{equation}
\mu _{a}(\mathbf{x})=\left\{ 
\begin{array}{cc}
c=const.>0, & \text{inside the tested inclusion,} \\ 
0, & \text{outside the tested inclusion.}%
\end{array}%
\right.  \label{6.4}
\end{equation}%
We perform the numerical tests with a variety of the value of the parameter $%
c=5,10,15,20,30$, see below. Therefore, by (\ref{6.3}) and (\ref{6.4}) 
\begin{equation}
\text{inclusion/background contrast}=1+\frac{c}{5}.  \label{6.5}
\end{equation}%
From the Physics standpoint $\mu _{s}(\mathbf{x})=5$ in the domain $\Omega $
means that an average particle scatters every 1/5 of the unit. Hence, by (%
\ref{6.0}) and (\ref{6.3}) the maximal average number of scattering events
for a particle emitted from a point $\mathbf{x}_{\alpha }\in \Gamma _{d}$ \
and entering $\Omega $ is around $5$ before this particle leaves $\Omega $ 
\cite{HT1,HT2,HT3}. This might happen in optics before the true diffusion
occurs, e.g. in the case of the so-called \textquotedblleft snake photons" 
\cite{Das}.

In computational results below $c\neq const.$ inside of computed shapes of
letters. Hence, we set for computed inclusion/background contrasts we
replace (\ref{6.5}) with:%
\begin{equation}
\text{computed inclusion/background contrast}=1+\frac{1}{5}\max_{\overline{%
\Omega }}\mu _{a}(\mathbf{x}).  \label{6.50}
\end{equation}

To solve the forward problem (\ref{2.7}), (\ref{2.14}), we have solved
integral equation (\ref{2.29}) with the function $u_{0}$ taken from (\ref%
{2.27}). To do this, we have used the discrete form of (\ref{2.27}), (\ref%
{2.29}) and the trapezoidal rule. The discretization steps with respect to $%
x,y,\alpha $ where $h_{x}=h_{y}=h_{\alpha }=1/40.$ The discretized integral
equation (\ref{2.29}) was solved as a linear system using the Matlab
function '\TEXTsymbol{\backslash}'. Thus, the solution of this forward
problem has provided computationally simulated data for the inverse problem
to us.

To minimize the convexification functional $J_{\lambda }\left( W^{h}\right) $
in (\ref{4.7}), we have written in the finite differences form not only the $%
x-$derivative as in section 4 but the $y-$derivative as well. Also,
integrals with respect to $\alpha $ were written in the discrete form using
the trapezoidal rule. The mesh sizes were different from ones for the
forward problem. They were: 
\begin{equation}
h_{x}=h_{y}=h_{\alpha }=h=1/20.  \label{6.6}
\end{equation}%
Then we have minimized the resulting functional $J_{\lambda ,dis}\left(
W^{h}\right) $, 
\begin{equation}
J_{\lambda ,dis}\left( W^{h}\right) =\left\Vert \left(
D_{N}^{h}W_{y}^{h}+A^{h}W_{x}^{h}+F^{h}\left( W^{h}\left( \mathbf{x}%
^{h}\right) ,\mathbf{x}^{h}\right) \right) e^{\lambda y}\right\Vert
_{L_{N}^{2,h}\left( \Omega ^{h}\right) }^{2}  \label{6.7}
\end{equation}%
in its fully discrete form with respect to the values of the vector function 
$W^{h}$ at the grid points. Vector functions and matrices in (\ref{6.7}) are
full analogs of those in (\ref{4.7}) with the only difference that they are
fully discrete in the above sense, rather than `partially' discrete as in
sections 4,5. The same is true for the norm $\left\Vert \cdot \right\Vert
_{L_{N}^{2,h}\left( \Omega ^{h}\right) }.$

It follows from (\ref{6.0}) and (\ref{6.6}) that we had total $20\times
20\times N$ unknown parameters in our minimization procedure. To solve the
minimization problem, we have used the Matlab's built-in function \textbf{%
fminunc} with the quasi-newton algorithm. The iterations of the function 
\textbf{fminunc} were stopping when the following inequality occurred at the
iteration number $k$: 
\begin{equation*}
\left\vert J_{\lambda ,dis}\left( W_{k}^{h}\right) \right\vert <10^{-2}.
\end{equation*}

By (\ref{3.22}) our technique requires computations of first derivatives
with respect to $\alpha .$ Note, however, that (\ref{3.10}) and (\ref{3.19})
imply that it is not necessary to calculate the $\alpha -$derivative of the
boundary data, which is an advantage, since boundary data are noisy. The
derivatives $\partial _{\alpha }$ of functions $\Psi _{s}\left( \alpha
\right) $ and the function $K(\mathbf{x},\alpha ,\beta )$ were calculated
via finite differences.

We have introduced the random noise in the boundary data $g_{1}(\mathbf{x}%
,\alpha )$ in (\ref{3.4}) on the boundary $\partial \Omega $, 
\begin{equation}
g_{1}(\mathbf{x},\alpha )=g_{1}(\mathbf{x},\alpha )\left( 1+\sigma \zeta _{%
\mathbf{x}}\right) .  \label{6.70}
\end{equation}%
Here $\zeta _{\mathbf{x}}$ is the uniformly distributed random variable in
the interval $[0,1]$ depending on the point $\mathbf{x}\in \partial \Omega $
with $\sigma =0.03$ and $\sigma =0.05,$ which correspond respectively to $3\%
$ and $5\%$ noise level.

To solve the minimization problem, we need to provide the starting $%
W_{0}^{h}(\mathbf{x}^{h})$ for iterations. Due to the global convergence
property of our method, the vector function $W_{0}^{h}(\mathbf{x}%
^{h})=\left( w_{0,0}^{h}(\mathbf{x}^{h}),...,w_{N-1,0}^{h}(\mathbf{x}%
^{h})\right) ^{T}$ should not have any information about the exact solution $%
W^{h\ast }(\mathbf{x}^{h}).$ On the other hand, due to (\ref{4.1}), we
should have $W_{0}^{h}\left( \mathbf{x}^{h}\right) \mid _{\partial \Omega
^{h}}=P^{h}\left( \mathbf{x}^{h}\right) .$ Therefore, in all numerical tests
below we choose the starting point as the discrete version of the following
vector function: 
\begin{equation*}
w_{s,0}(x,y)=\frac{1}{2}\left( \frac{(A-x)}{2A}w_{s}(-A,y)+\frac{(x+A)}{2A}%
w_{s}(A,y)\right) +
\end{equation*}

\begin{equation}
+\frac{1}{2}\left( \frac{(b-y)}{b-a}w_{s}(x,a)+\frac{(y-a)}{b-a}%
w_{s}(x,b)\right) ,s=0,...,N-1.  \label{6.8}
\end{equation}
Expression (\ref{6.8}) represents the average of linear interpolations
inside of the square $\Omega $ with respect to $x-$direction and $y-$%
direction of the boundary condition for $w_{s}\left( x,y\right) $ $.$

\subsection{Numerical results}

\label{sec:6.2}

Recall that we reconstruct the coefficient $a\left( \mathbf{x}\right) =\mu
_{s}\left( \mathbf{x}\right) +\mu _{a}\left( \mathbf{x}\right) ,$ where
functions $\mu _{s}\left( \mathbf{x}\right) $ and $\mu _{a}\left( \mathbf{x}%
\right) $ are given in (\ref{6.3}) and (\ref{6.4}) respectively. Our results
for Tests 1-4 are for noiseless data and the results for Test 5 are for
noisy data as in (\ref{6.70}).

To demonstrate a good performance of our technique, we intentionally test it
for rather complicated shapes of inclusions, which are non convex and have
voids. More precisely, our inclusions are letters `A', `$\Omega $' and two
letters jointly `SZ'. `SZ' stands for `Shenzhen', the city where the campus
of the Southern University of Science and Technology, the work place of the
third and fourth authors, is located. Thus, in our tests the coefficients $%
\mu _{a}(\mathbf{x})$ in (\ref{6.4}) have shapes of those letters located
inside of the $1\times 1$ square $\Omega $ defined in (\ref{6.0}).

\textbf{Test 1}. We test the letter `A' with $c=5$ in (\ref{6.4}). This is
our reference case. More precisely, we use this test to figure out optimal
values of parameters $N$ and $\lambda .$ As soon as optimal parameters are
selected, we use them then for all other tests.

First, we select $N$. To do this, we solve the forward problem (\ref{2.7})%
\emph{, }(\ref{2.14}) for the case when the functions $\mu _{s}\left( 
\mathbf{x}\right) $ and (\ref{6.3})\emph{\ }and (\ref{6.4}) respectively and
in $c=5.$ Hence, by (\ref{6.5}) the inclusion/background contrast is 2:1 in
this case. Next, we calculate norms $\left\Vert w_{s}\left( \mathbf{x}%
\right) \right\Vert _{L_{2}\left( \Omega \right) }$ and compare them. We
have observed that the $L_{2}\left( \Omega \right) -$norm of the function $%
w_{s}(\mathbf{x})$ decreases very rapidly when the number $s$ is growing.
More precisely, we have obtained that 
\begin{equation}
\frac{\dsum\limits_{s=3}^{11}\left\Vert w_{s}\left( \mathbf{x}\right)
\right\Vert _{L_{2}\left( \Omega \right) }}{\dsum\limits_{s=0}^{11}\left%
\Vert w_{s}\left( \mathbf{x}\right) \right\Vert _{L_{2}\left( \Omega \right)
}}=0.0084,  \label{6.10}
\end{equation}%
which means less than 1\%. The values of norms $\left\Vert w_{s}\left( 
\mathbf{x}\right) \right\Vert _{L_{2}\left( \Omega \right) }$ for $s=0,...,11
$ are displayed in Table \ref{table_u_norm_basis} and Figure \ref%
{plot_u_norm_basis}. One can observe that starting from $s=3$, these norms
are much less than those for $s=0,1,2.$ We conclude therefore, that we
should take in our tests%
\begin{equation}
N=3.  \label{6.9}
\end{equation}

\begin{table}[htbp]
\caption{The $L_{2}\left( \Omega \right) -$norms of functions $w_{s}\left( 
\mathbf{x}\right)$, $s=0,1,...,11$ for the reference Test 1 with $c=5$ in (%
\protect\ref{6.4}).}
\label{table_u_norm_basis}\centering
\begin{tabular}{c|c|c|c|c|c|c|}
\hline
$s$ & 0 & 1 & 2 & 3 & 4 & 5 \\ \hline
$\left\| w_{s}(\mathbf{x}) \right\|_{L_{2}} $ & 5.7122 & 1.6383 & 0.1630 & 
0.0118 & 0.0091 & 0.0077 \\ \hline
$s$ & 6 & 7 & 8 & 9 & 10 & 11 \\ \hline
$\left\| w_{s}(\mathbf{x}) \right\|_{L_{2}} $ & 0.0067 & 0.0061 & 0.0055 & 
0.0057 & 0.0058 & 0.0054 \\ \hline
\end{tabular}%
\end{table}

\begin{figure}[tbph]
\centering
\includegraphics[width = 3.5in]{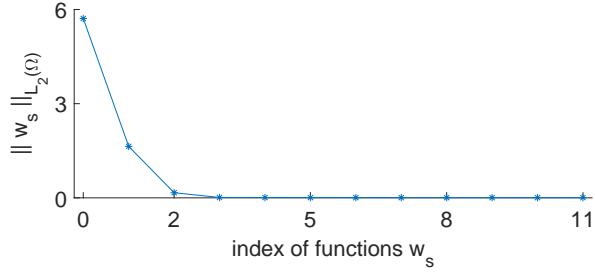}
\caption{The decrease with respect to $s$ of the $L_{2}\left( \Omega \right)
-$norms of functions $w_{s}\left( \mathbf{x}\right)$, $s=0,1,...,11$ for the
reference Test 1 with $c=5$ in (\protect\ref{6.4}).}
\label{plot_u_norm_basis}
\end{figure}

\begin{figure}[tbph]
\centering
\includegraphics[width = 5in]{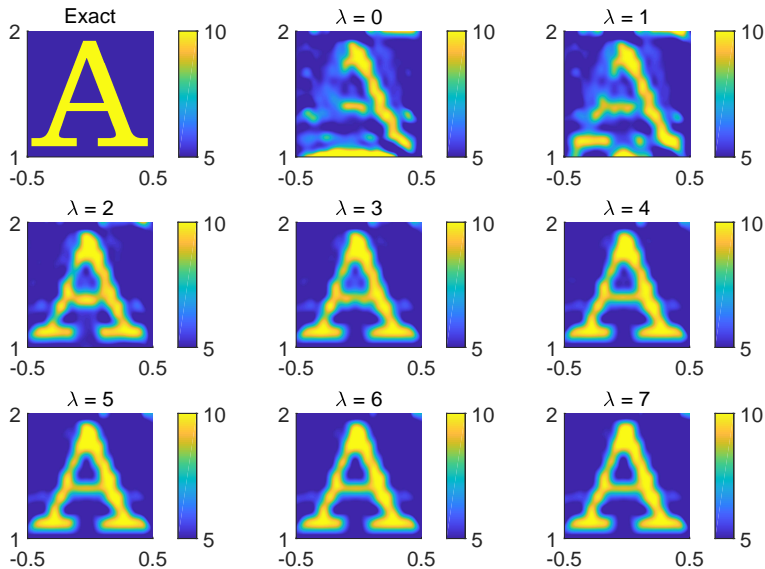}
\caption{The reconstructed coefficient $a\left( \mathbf{x}\right) $, where
the function $\protect\mu_{a}\left( \mathbf{x}\right) $ is given in (\protect
\ref{6.4}) with $c=5$ inside of the letter `A'. The goal of here is to test
different values of the parameter $\protect\lambda =0,1,2,3,4,5,6,7$ for $N=3
$ as in (\protect\ref{6.9}). The value of $\protect\lambda $ can be seen on
the top side of each square. The images have a low quality for $\protect%
\lambda =0, 1.$ Then the quality is improved and is stabilized at $\protect%
\lambda =5.$ Thus, we select $\protect\lambda =5$ as an optimal value of
this parameter for all follow up tests.}
\label{plot_NBasis3_diff_Lambda}
\end{figure}

Next, given the optimal value of $N=3$ in (\ref{6.9}), we select the optimal
value of the parameter $\lambda $ of the Carleman Weight Function $%
e^{2\lambda y}$ in (\ref{6.7}). To do this, we test the same letter `A' with 
$c=5$ inside of it for values of the parameter $\lambda =0,1,2,3,4,5,6,7.$
Our numerical results are presented on Figure \ref{plot_NBasis3_diff_Lambda}%
. We observe that the images have a low quality for $\lambda =0,1.$ Then the
quality is improved and is stabilized at $\lambda =5.$ Thus, we treat $%
\lambda =5$ as the optimal value of this parameter. We use this value in all
subsequent tests. 
The value $\lambda =5$ tells us that even though our
theorems 4.1-4.5 require sufficiently large values of the parameter $\lambda
,$ the computational practice shows that a reasonable value of $\lambda $
can be chosen. The same observation was made in all previous works on the
convexification of this research group, see items 2 and 3 in Remarks 4.1 in
the end of section 4. 

\begin{figure}[tbph]
\centering
\includegraphics[width = 5in]{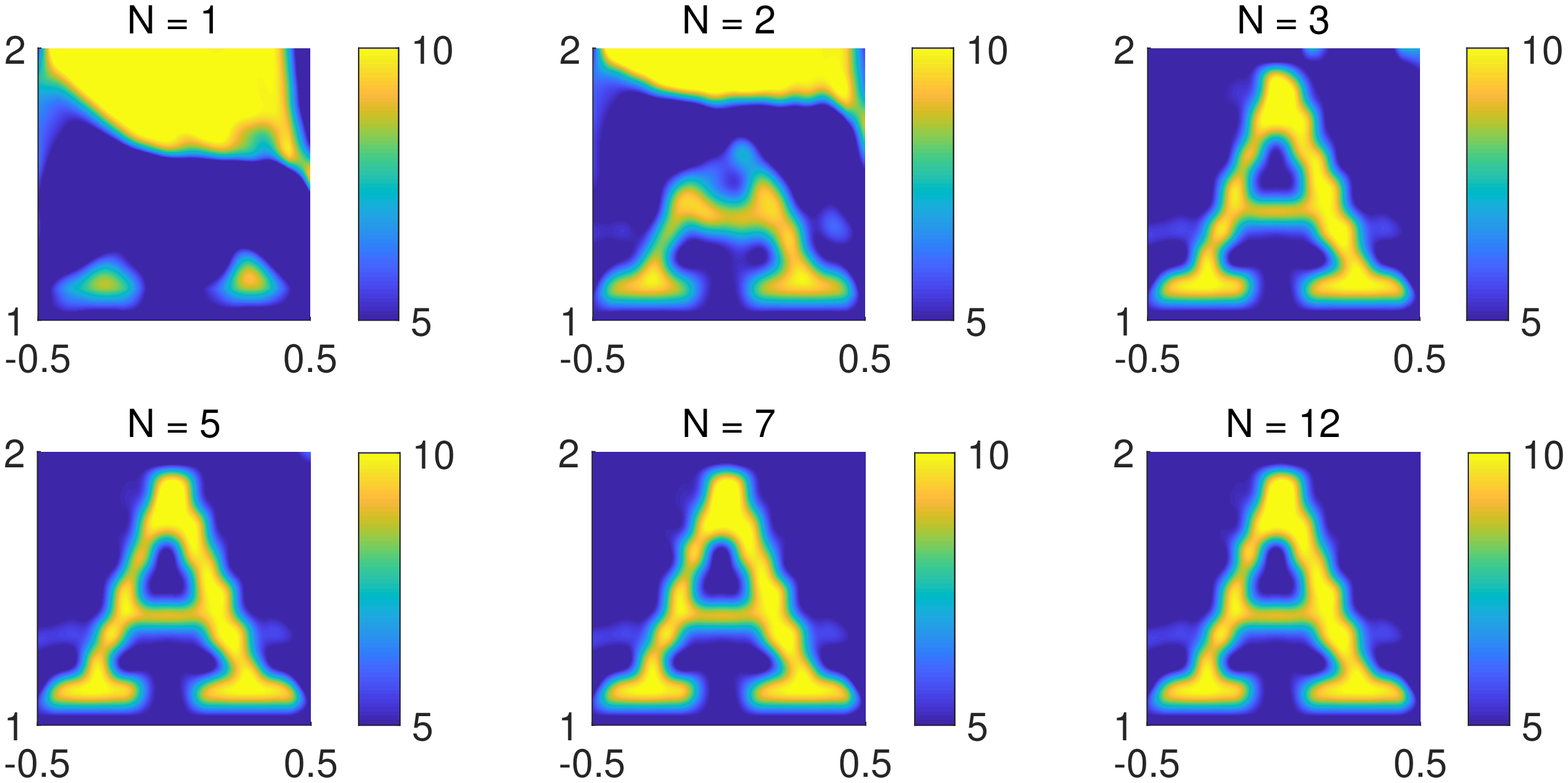}
\caption{The reconstructed coefficient $a\left( \mathbf{x}\right) $, where
the function $\protect\mu _{a}\left( \mathbf{x}\right) $ is given in (%
\protect\ref{6.4}) with $c=5$ inside of the letter `A'. We took the optimal
value of the parameter $\protect\lambda =5$ (see Figure \protect\ref%
{plot_NBasis3_diff_Lambda}) and have tested different values of the
parameter $N=1,2,3,5,7,12$. A low quality can be observed for $N=1,2$. The
reconstructions are basically the same for $N=3,5,7,12$. However, the
computational cost increases very rapidly with the increase of $N$, which is
explained by (\protect\ref{6.10}), Table \protect\ref{table_u_norm_basis}
and Figure \protect\ref{plot_u_norm_basis}. We conclude, therefore, that to
balance between the reconstruction accuracy and the computational cost, we
should use $N=3$ as in (\protect\ref{6.9}). Thus, we use below $\protect%
\lambda =5$ and $N=3$.}
\label{plot_lambda5_diff_NBasis}
\end{figure}

Now we want to demonstrate numerically again that $N=3$ is indeed a good
choice of $N$ for our optimal value of $\lambda =5.$ Taking $\lambda =5,$ we
test the same letter `A' as above with $c=5$ in it, but for $N=1,2,3,5,7,12.$
The results are displayed in Figure \ref{plot_lambda5_diff_NBasis}. One can
observe that reconstructions have a low quality for $N=1,2$. $\ $Next, the
reconstructions are basically the same for $N=3,5,7,12.$ However, the
computational cost increases very rapidly with the increase of $N$. Thus,
using also Table \ref{table_u_norm_basis}, Figure \ref{plot_u_norm_basis}
and (\ref{6.9}), we conclude that to balance between the reconstruction
accuracy and the computational cost, we should use $N=3$. Thus, in all
subsequent computations we use 
\begin{equation}
N=3,\lambda =5.  \label{6.11}
\end{equation}%

\textbf{Test 2}. We test the reconstruction of the coefficient $a(\mathbf{x})
$ with the shape of the letter `A' where the function $\mu _{a}\left( 
\mathbf{x}\right) $ is given in (\ref{6.4}) with different values of the
parameter $c=15,20,30$ inside of the letter `A'. Thus, by (\ref{6.50}) the
inclusion/background contrasts now are respectively $4:1$, $5:1$ and $6:1$.
Our computational results for this test are displayed on Figure \ref%
{plot_A10_A15_A20_A30}. One can observe that the quality of these images is
good for all four cases, although it slightly deteriorates for $c=20$ and $%
c=30.$ The computed inclusion/background contrast is accurate, see (\ref%
{6.50}). 
\begin{figure}[tbph]
	\centering
	\includegraphics[width = 3in]{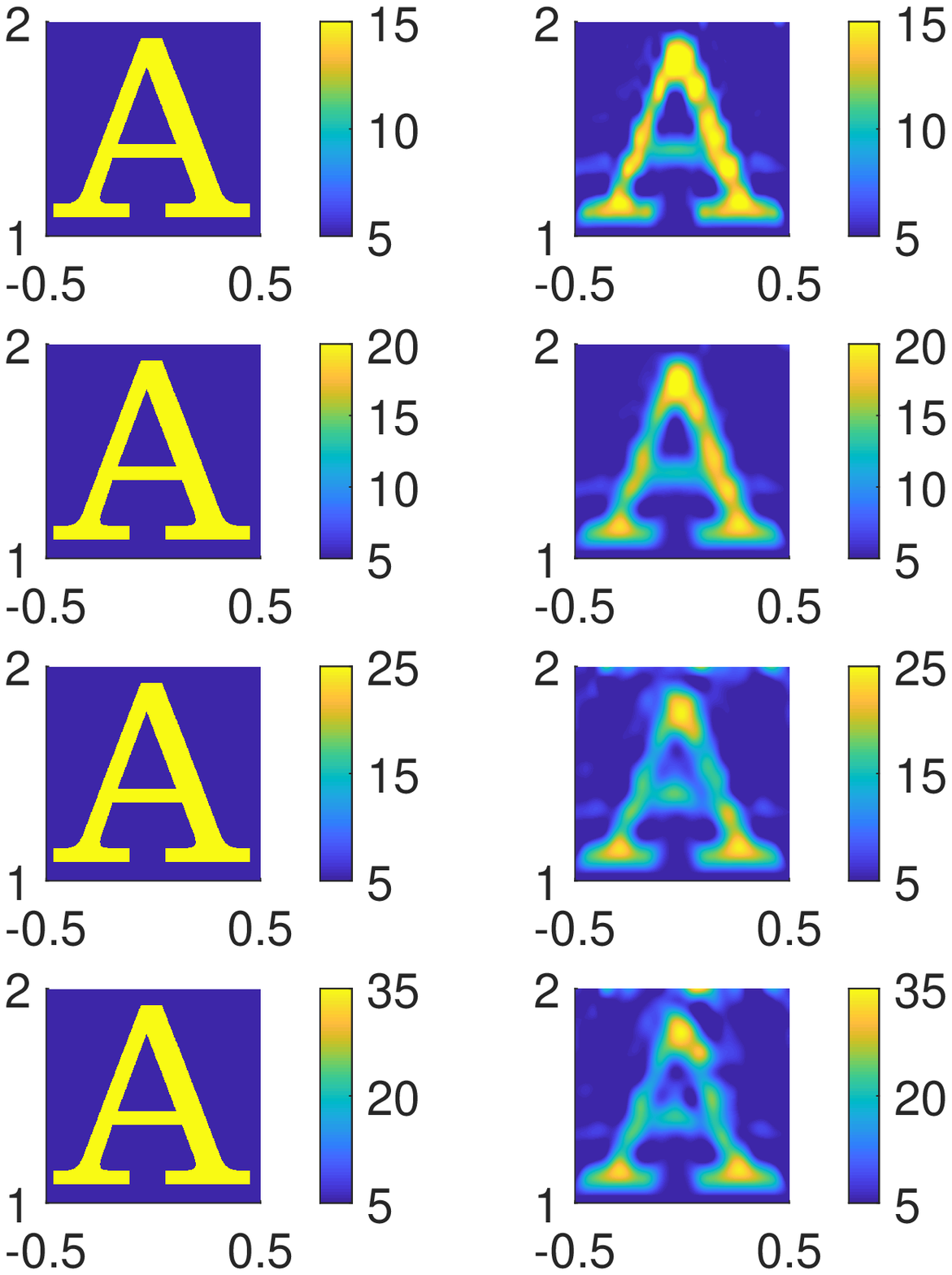}  
	\caption{Test 2. Exact (left) and reconstructed (right) coefficient $a(\mathbf{x})$ for $c=10,15,20,30$ inside of the letter `A' as in (\ref{6.4}) for $N=3,\lambda =5$, see (\ref{6.11}). 
	Thus, by (\ref{6.50}) the inclusion/background contrasts now are respectively $4:1$, $5:1$ and $6:1$. 
	The image quality remains basically the same for all these values of the parameter $c$, although some deterioration of this quality can be observed for $c=20$ and $c=30$. 
	The computed inclusion/background contrasts (\ref{6.50}) are accurate.}
	\label{plot_A10_A15_A20_A30}
\end{figure}

\begin{figure}[tbph]
\centering
\includegraphics[width = 3.5in]{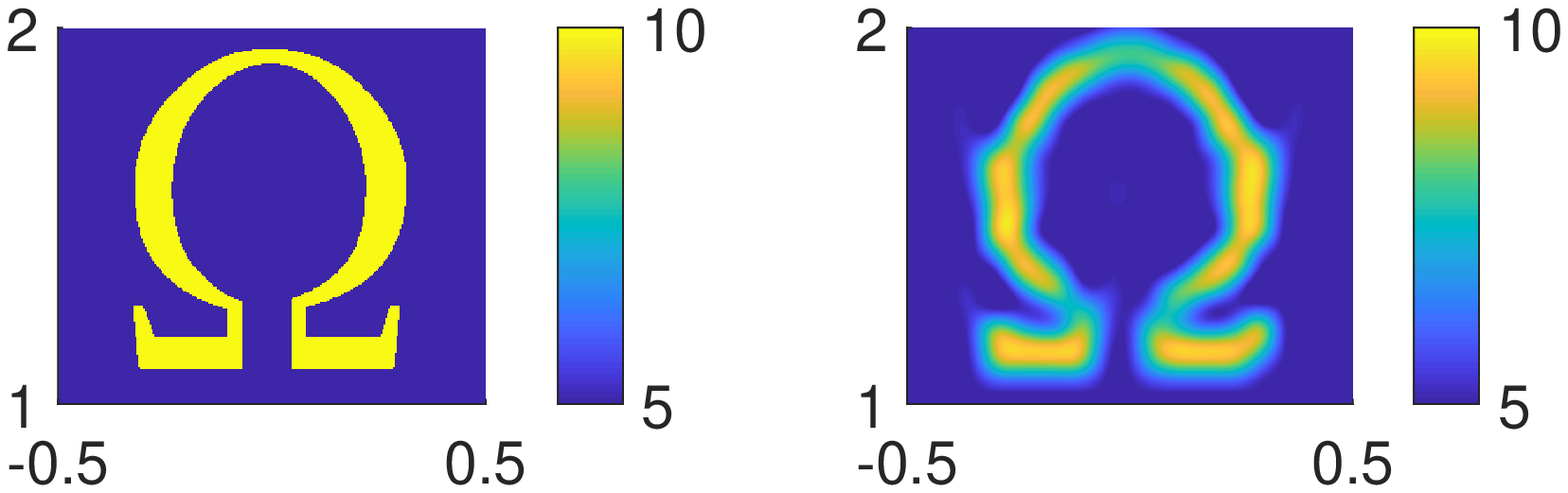}
\caption{Test 3. Exact (left) and reconstructed (right) coefficient $a(%
\mathbf{x})$ for the case when the function $\protect\mu _{a}\left( \mathbf{x%
}\right) $ is given in (\protect\ref{6.4}) with $c=5$ inside of the letter `$%
\Omega $'. The reconstruction is accurate.}
\label{plot_Omega}
\end{figure}

\textbf{Test 3}. We test the reconstruction of the coefficient $a(\mathbf{x}%
) $ with the shape of the letter `$\Omega $' where the function $\mu
_{a}\left( \mathbf{x}\right) $ is given in (\ref{6.4}) with $c=5$ inside of
the letter `$\Omega $'. Results are presented on Figure \ref{plot_Omega}. We
again observe an accurate reconstruction.

\begin{figure}[tbph]
\centering
\includegraphics[width = 3.5in]{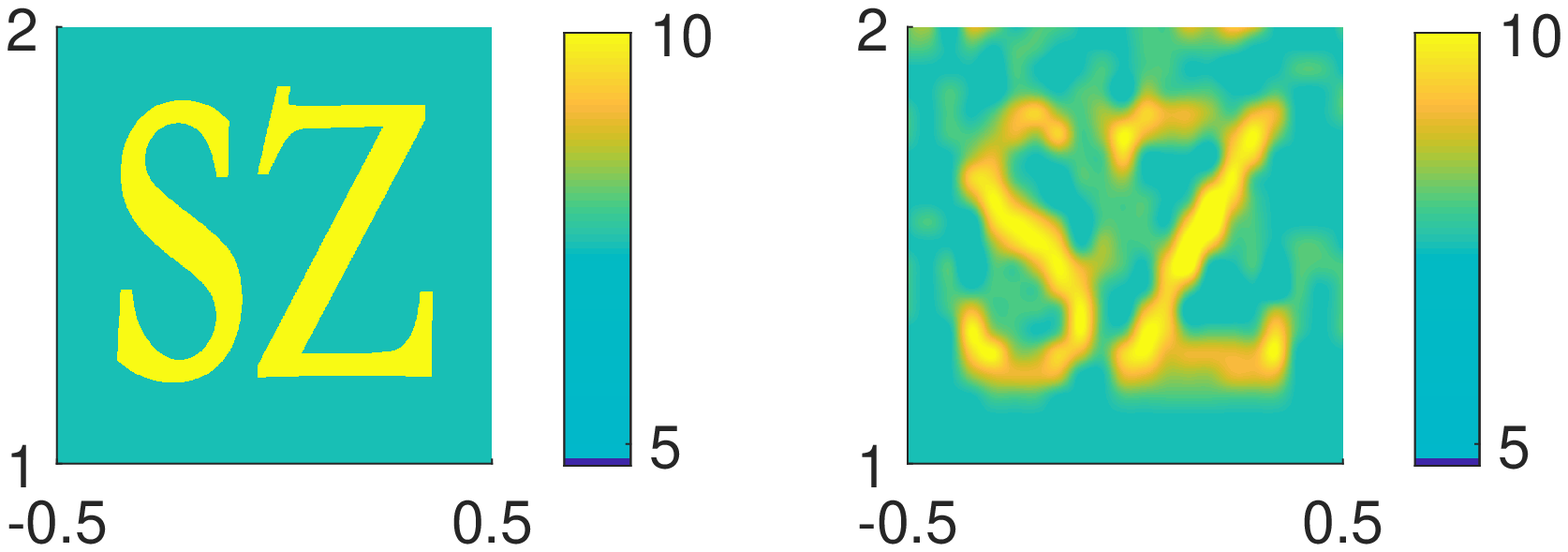}
\caption{Test 4. Exact (left) and reconstructed (right) coefficient $a(%
\mathbf{x})$ for the case when the function $\protect\mu _{a}\left( \mathbf{x%
}\right) $ is given in (\protect\ref{6.4}) with $c=5$ with the shape of two
letters `SZ'. In (\protect\ref{6.4}) $c=5$ inside of each of these two
letters and $\protect\mu _{a}\left( \mathbf{x}\right) =0$ outside of each of
these two letters. Here $N=3,\protect\lambda =5$ as in (\protect\ref{6.11}).
The image quality is lower than one for the case of the single letter $%
\Omega $ on Figure \protect\ref{plot_Omega}. 
Nevertheless, the quality is still good and the computed inclusion/background contrasts (\ref{6.50}) are accurate in both letters.}
\label{plot_SZ}
\end{figure}

\textbf{Test 4}. We test the reconstruction of the coefficient $a(\mathbf{x})
$ with the shape of two letters `SZ' where the function $\mu _{a}\left( 
\mathbf{x}\right) $ is given in (\ref{6.4}) with $c=5$ inside of each of
these two letters and $\mu _{a}\left( \mathbf{x}\right) =0$ outside of each
of these two letters. In this test, $N=3,\lambda =5$ as in (\ref{6.11}).
Results are presented on Figure \ref{plot_SZ}. The image quality is lower
than one for the case of the single letter `$\Omega $' on Figure \ref%
{plot_Omega}. Nevertheless, the quality is still good and the computed
inclusion/background contrasts (\ref{6.50}) are accurate in both letters.

\begin{figure}[tbph]
\centering
\includegraphics[width = 3.5in]{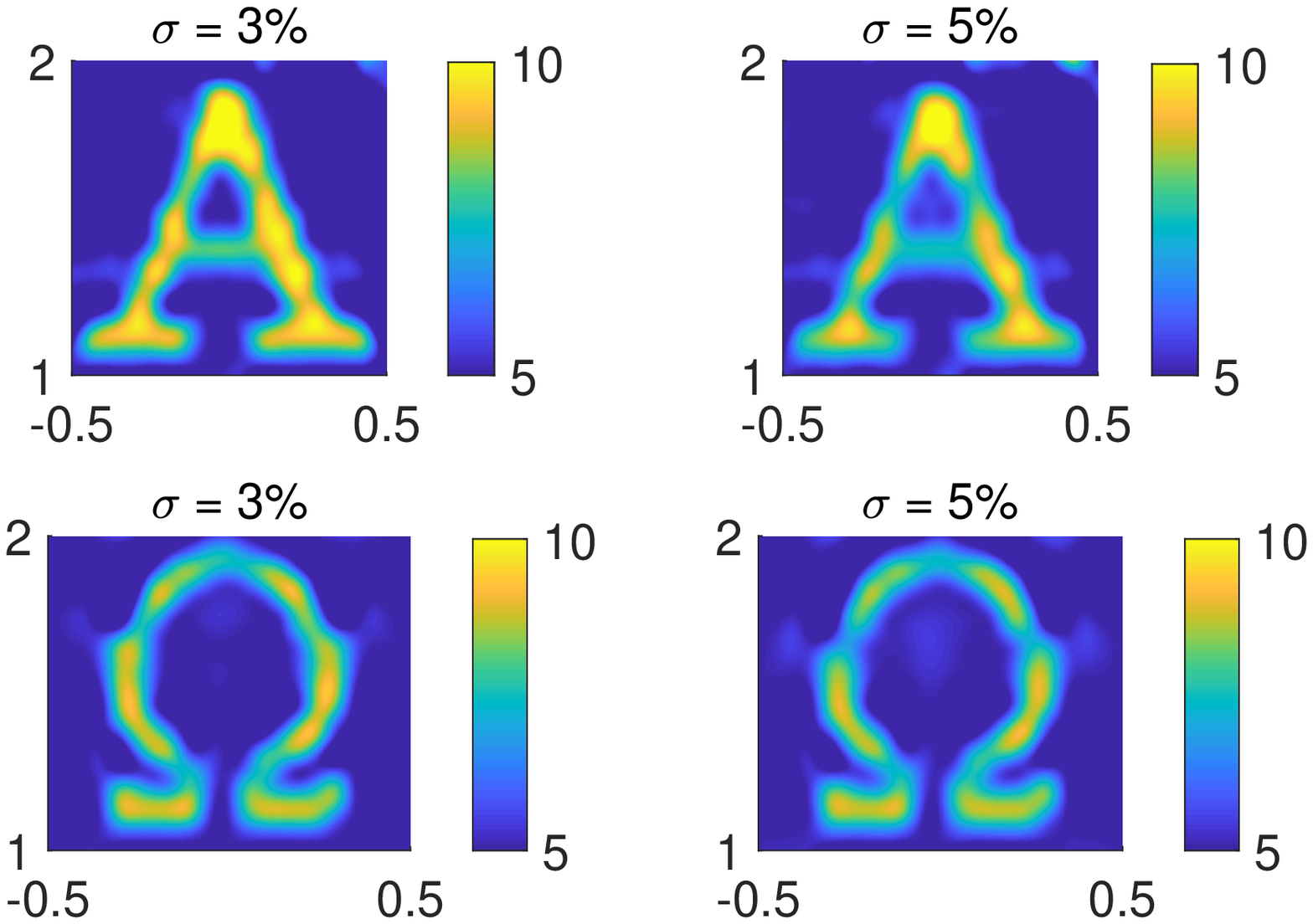}
\caption{Reconstructed coefficient $a(\mathbf{x})$ with the shape of letters
`A' (top) and `$\Omega $' (bottom) with $c=5$ from noise polluted
observation data as in (\protect\ref{6.70}) with $\protect\sigma =0.03$
(left) and $\protect\sigma =0.05$ (right), i.e. with 3\% (left) and 5\%
(right) noise level. Here $N=3$ and $\protect\lambda =5$ as in (\protect\ref%
{6.11}).
One can observe accurate reconstructions in all four cases. In particular, the inclusion/background contrasts (\ref{6.50}) are reconstructed accurately.}
\label{plot_AddNoise}
\end{figure}

\textbf{Test 5}. In this test we use noisy data as in (\ref{6.70}) with $%
\sigma =0.03$ and $\sigma =0.05,$ i.e. with 3\% and 5\% noise level. We test
the reconstruction of the coefficient $a(\mathbf{x})$ with the shape of
either the letter `A' or the letter `$\Omega $', where the function $\mu
_{a}\left( \mathbf{x}\right) $ is given in (\ref{6.4}) with $c=5$ inside of
each of these two letters. Again, $N=3$, $\lambda =5$ as in (\ref{6.11}).
The results are shown in Figure \ref{plot_AddNoise}. One can observe
accurate reconstructions in all four cases. In particular, the
inclusion/background contrasts (\ref{6.50}) are reconstructed accurately.

\providecommand{\bysame}{\leavevmode\hbox to3em{\hrulefill}\thinspace} %
\providecommand{\MR}{\relax\ifhmode\unskip\space\fi MR } 
\providecommand{\MRhref}[2]{  \href{http://www.ams.org/mathscinet-getitem?mr=#1}{#2}
} \providecommand{\href}[2]{#2}


\end{document}